\documentclass[11pt,reqno]{amsart}                                                                                                                                                                                                     

\usepackage{                                                                                                                                                                                                                           
amssymb,                                                                                                                                                                                                                               
amsaddr                                                                                                                                                                                                                                
}                                                                                                                                                                                                                                      
\newcommand{\no}[1]{#1}                                                                                                                                                                                                                

\renewcommand{\no}[1]{}  \newcommand{\upDelta}{\Delta} 

\no{\usepackage{times}\usepackage[                                                                                                                                                                                                     
subscriptcorrection, slantedGreek,                                                                                                                                                                                                     
nofontinfo]{mtpro}                                                                                                                                                                                                                     
\renewcommand{\Delta}{\upDelta}                                                                                                                                                                                                        
}                                                                                                                                                                                                                                      

\date{\today}                                                                                                                                                                                                                          

\usepackage{algpseudocode}                                                                                                                                                                                                             
\usepackage{algorithm}                                                                                                                                                                                                                 
\usepackage{color}
\usepackage{graphicx}
\usepackage{epstopdf}
\usepackage{subfigure}
\usepackage[normalem]{ulem}

\newtheorem{lemma}{Lemma}                                                                                                                                                                                                              
\newtheorem{definition}{Definition}                                                                                                                                                                                                    
\newtheorem{corollary}{Corollary}

\theoremstyle{remark}

\DeclareMathOperator{\sinc}{sinc}

\newcommand{\be}[1]{\begin{equation}\label{#1}}                                                                                                                                                                                        
\newcommand{\ee}{\end{equation}}                                                                                                                                                                                                       

\newcommand{\diag}{\mathop{\mathrm{diag}}}

\title{A continuous adjoint for photo-acoustic tomography of the brain}

\author[A. Javaherian]{Ashkan Javaherian and Sean Holman}                                                                                                                                                                                              

\address[]{School of Mathematics \\ University of Manchester  \\ Manchester, M19 7PL UK \\ (e-mail: ashkan.javaherian@postgrad.manchester.ac.uk)}                                                                                      
\date{}

\begin{document}

\maketitle
\begin{abstract}
We present an optimization framework for photo-acoustic tomography of brain based on a system of coupled equations that describe the propagation of sound waves in linear isotropic inhomogeneous and lossy elastic media with the absorption and physical dispersion following a frequency power law using fractional Laplacian operators. The adjoint of the associated continuous forward operator is derived, and a numerical framework for computing this adjoint based on a k-space pseudo-spectral method is presented. We analytically show that the derived continuous adjoint matches the adjoint of an associated discretised operator. We include this adjoint in a first-order positivity constrained optimization algorithm that is regularized by total variation minimization, and show that the iterates monotonically converge to a minimizer of an objective function, even in the presence of some error in estimating the physical parameters of the medium.
\end{abstract}
\section{Introduction}\label{Intro}                                                                                                                                                                                                    

Quantitative Photo-acoustic Tomography (QPAT) is a hybrid imaging modality which simultaneously takes advantage of the rich contrast attributed to optical imaging and the high spatial resolution brought up by ultrasound. In this technique, short pulses of near-infrared light are used to irradiate tissue. The energy from these pulses is absorbed as a function of the optical absorption map of the tissue. This generates local increases in pressure which propagate outwards as photo-acoustic (PA) waves, and are then measured by broadband detectors placed at the surface. The inverse problem of QPAT is to reconstruct the spatially varying optical absorption coefficient from the recorded PA signals. This involves two inverse problems, namely acoustic and optical \cite{Rosenthal}. These two inverse problems can be solved distinctly \cite{Huang,Arridge,Javaherian}, or alternatively as a direct hybrid problem \cite{Haltmeier}. In this work we consider only the acoustic portion of the inverse problem which we simply call Photo-acoustic Tomography (PAT).

Considering the acoustic inverse problem, Time reversal (TR) is a comprehensive inversion approach for PAT since it can be used for media with heterogeneous acoustic properties
and arbitrary detection geometries \cite{Hristova-a,Hristova-b,Treeby}. However, this method is based on a continuous domain with idealized conditions such as a closed detection surface or exactly known medium's properties \cite{Hristova-a}, which do not hold in real cases. Problems such as finite sampling, a limited accessible angle for detection surface, errors in estimation of medium's properties, or errors in data measurement make the acoustic inverse problem ill-posed \cite{Huang}. In these cases, model-based iterative methods are often used, e.g., TR-based iterative techniques \cite{Qian}, or optimization algorithms \cite{Huang,Javaherian}. The optimization approaches are often based on computation of
the gradient of an objective function in terms of a forward model and the corresponding adjoint model.                                                                                                                                       

Because of the dependance of shape, spectrum and amplitude of PA signals on physical properties of tissue media, it will be advantageous if the image reconstruction in PAT is enriched by tissue-realistic models that account for the absorption behaviours evident in tissues \cite{Treeby-a,Haltmeier-b,Scherzer}. Among model-based iterative approaches for absorbing media, the adjoint was computed by a “discretize-then-adjoint” method in \cite{Huang}, or by an “adjoint-then-discretize” method in \cite{Javaherian}.

It is well-known that modelling the propagation of sound waves can be considerably expedited compared to Finite difference time-domain (FDTD) methods by using Pseudo-spectral time-domain (PSTD) methods. Applying these techniques, the spatial gradients are computed in frequency domain, while the temporal gradients are computed using finite difference methods, similar to FDTD techniques. The efficiency of PSTD methods is because of a fast computation of the spatial gradients using Fast Fourier Transforms (FFTs), as well as a dramatic relaxation in the mesh requirement and time step \cite{Cox,Firouzi}.

In PAT, the compartmentalised distribution of light absorbing molecules composing tissues induces step-like discontinuities in the generated pressure field. As a result, the generated PA waves are considerably more broadband than ultrasonic waves \cite{Treeby,Javaherian}. Furthermore, the absorption of sound waves in many media such as tissues has been experimentally shown to follow a frequency power law with a non-integer power, which can be described by fractional derivatives \cite{Treeby-a,Treeby-f}. Classical attenuation models used the fractional time derivatives, which are non-local in time, and thus require storing the time history of field variables \cite{Mainard}. It has been established that the fractional time derivatives can be replaced by fractional space derivatives, which are nonlocal in space rather than time, and are thus more memory efficient \cite{Treeby-a,Treeby-f}. This is done using the dispersion relation for lossless wave equation. The cost of this method is that the spatially non-local operators violate causality \cite{Kowar}.

PAT has shown its potential for characterization of the vasculature in small animals or within a few mm of the skin's surface in humans \cite{Zhang}. Furthermore, PAT has been utilized successfully for transcranial brain imaging in small animals \cite{Xu-i,Li}. In these cases, the effect of the skull on the propagation of PA waves is neglected because of the low thickness of the skull ($\approx 1\  \text{mm}$), and thus the image reconstruction is done based on scalar acoustic
wave equations \cite{Xu-i,Li}.

To account for aberration of PA signals because of the heterogeneous properties of the skull, a subject-specific imaging model was proposed, where the inhomogeneity of the skull is taken into account using adjunct information about the skull anatomy and composition \cite{Huang-a}. This information must be obtained from x-ray computed tomography image data, or some other imaging modalities.

The application of PAT in transcranial brain imaging of humans is very limited since PA signals are aberrated to a high degree by absorption, scattering and compressional-to-shear mode conversion effects due to the high thickness of the skull
(4mm-7mm). Recently, a numerical framework for image reconstruction in transcranial PAT was proposed, where the forward problem describes the wave propagation in a linear isotropic, heterogeneous and lossy elastic medium, and the corresponding adjoint model is obtained by an explicit reversal of the computational steps of the forward solver, i.e.,
\textit{discretize-then-adjoint} method \cite{Mitsuhashi}. The adjoint derived by this technique is not the adjoint of the continuous forward model, but of the particular numerical scheme. This forward and adjoint pair was discretised using the finite-difference time domain (FDTD) method, and the attenuation effects were described by a diffusive model, which ignores the dependency of the wavefield attenuation on frequency \cite{Mitsuhashi}. Consider that any changes in the forward model, e.g., using tissue-realistic absorption models or high-performance solvers, require the reformulation of the algebraic adjoint.

In elastic solids, compressional and shear waves propagate at different speeds. As a result, using the dispersion relation for describing fractional space derivatives requires splitting
the field variables into compressional and shear parts \cite{Treeby-f}. This is done using a dyadic wave number tensor in the frequency domain \cite{Firouzi,Treeby-f}. Additionally, by splitting the fields, the numerical dispersion errors accumulated by the time integrations can also be minimized via applying the k-space correction to the spatial gradients, which allows larger time steps without loss of stabilty or accuracy in heterogeneous media \cite{Tabei,Firouzi}.

\textit{Contribution.} We consider a forward map in the PAT problem in which a system of coupled first-order equations describes the propagation of PA waves in linear isotropic, heterogeneous and lossy elastic media, where the abosrption and physical dispersion follow a frequency-power law. We derive the adjoint of the PAT forward map in this context. This adjoint, referred to here as the \textit{analytic viscoelastic adjoint}, is derived on a continuous domain, and is in the form of a system of partial differential equations. One of the advantages of the analytic adjoint over the algebraic adjoint derived in \cite{Mitsuhashi} is that this is agnostic to the numerical scheme used for solving the equations. Another advantage is that by setting viscosity coefficients to zero in the derived analytic adjoint, the general form of the adjoint model for lossless media is derived. This can be used as a basic model, when other existing attenuation models are considered, e.g. \cite{Kowar,Kowar-b}. We shall analytically show that a numerical computation of the derived \textit{analytic viscoelastic adjoint} using the k-space pseudo-spectral method matches the algebraic adjoint of the associated forward model. The derived analytic adjoint is numerically validated using the adjoint test, and then the forward and adjoint pair is included in a positivity constrained and total-variation regularized solver based on the Iterative Shrinkage Thresholding algorithm (ISTA) for image reconstruction in 2D and 3D scenarios \cite{Beck}.

\section{Background}        \label{s2}                                                                                                                                                                                                           

The relation between stress tensor $\sigma$ and strain tensor $\epsilon$ in an isotropic lossless elastic medium is described using the Einstein summation notation in the form
\begin{align}
\sigma_{ij}  = \lambda\delta_{ij} \epsilon_{ll}+2\mu \epsilon_{ij}
\end{align}
in Cartesian coordinates. Here, $\mu$ and $\lambda$ are the Lam\'e elastic parameters, and are related to the shear and compressional wave speeds, $c_s$ and $c_p$ respectively, by the equations
\begin{align}
\mu= \rho c_s^2, \hspace{0.5cm} \lambda= \rho c_p^2-2 \mu,
\end{align}
where $\rho$ denotes the medium's mass density.
The strain tensor is a function of the particle displacement vector ($u$) in the form
\begin{align} \label{str-dis}                                                                                                                                                                                                          \epsilon_{ij}=\frac{1}{2}\left(\frac{\partial u_i}{\partial x_j}+\frac{\partial u_j}{\partial x_i}\right),                                                                                                                             \end{align}
where $x$ stands for the position.

For an isotropic viscoelastic medium the stress-strain relationship can be described by a classical variant of the so-called Kelvin-Voigt model, which accounts for an acoustic absorption proportional to $\omega^2$ and no dispersion in the low-frequency limit, where $\omega$ denotes the temporal frequency \cite{Treeby-f}. This model is in the form
\begin{align} \label{Kelvin}                                                                                                                                                                                                           \sigma_{ij}  &=\lambda\delta_{ij} \epsilon_{ll}+2\mu \epsilon_{ij}+\chi \delta_{ij}\frac{\partial}{\partial t} \epsilon_{ll}+2\eta \frac{\partial}{\partial t}\epsilon_{ij},
\end{align}                                                                                                                                                                                                                            where $\chi$ and $\eta$ denote the compressional and shear viscosity coefficients. Setting $\chi,\eta=0$ gives
the stress-strain relation for lossless media. Plugging \eqref{str-dis} into \eqref{Kelvin}, together with $v=du/dt$ with $v$ denoting the particle velocity vector, gives                                                                                                                                                           
\begin{align} \label{dis-Kelvin}                                                                                                                                                                                                       \frac{\partial \sigma_{ij}}{\partial t}=\lambda\delta_{ij}\frac{\partial v_l}{\partial x_l} +\mu \left(\frac{\partial v_i}{\partial x_j}+\frac{\partial v_j}{\partial x_i}\right) +\chi \delta_{ij} \frac{\partial^2 v_l} {\partial
x_l \partial t}+ \eta \left( \frac{\partial^2 v_i} {\partial x_j \partial t}+ \frac{\partial^2 v_j} {\partial x_i \partial t}\right).
\end{align}                                                                                                                                                                                                                            
However, as discussed in section \ref{Intro}, experimental studies have shown that attenuation in many materials of interests, including tissue media such as bone, is proportional to $\omega^y$ with $y$ a non-integer between 0 and 2 \cite{Szabo}. Because of the broadband nature of PA signals, as well as the high level of the attenuation in the skull, this behaviour cannot be neglected. To account for the non-integer power law dependence, the integer temporal derivatives in equations \eqref{Kelvin} and \eqref{dis-Kelvin} can be replaced by fractional time derivatives \cite{Holm}. For an isotropic medium, this gives the fractional Kelvin-Voigt model in the form
\begin{align} \label{frKelvin}                                                                                                                                                                                                         \sigma_{ij}  &=\lambda\delta_{ij} \epsilon_{ll}+2\mu \epsilon_{ij}+\chi \delta_{ij} \frac{\partial^{y-1}}{\partial t^{y-1}} \epsilon_{ll} +  2\eta \frac{\partial^{y-1}}{\partial t^{y-1}}\epsilon_{ij},
\end{align}                                                                                                                                                                                                                            
where                                                                                                                                                                                                                                  
\begin{align}                                                                                                                                                                                                                          \eta=-\frac{2\rho c_s^3}{\cos{(\pi y/2)}} \alpha_{0,s}, \hspace{0.5cm} \chi= -\frac{2\rho c_p^3}{\cos{(\pi y/2)}} \alpha_{0,p}- 2\eta,
\end{align}                                                                                                                                                                                                                            
with $\alpha_{0,s}$ and $\alpha_{0,p}$, respectively the attenuation coefficients pertaining to shear and compressional waves in $\text{Np} (\text{rad}/\text{s})^{-y} \text{m}^{-1}$ \cite{Treeby-f}.

The temporal fractional derivatives in equation \eqref{frKelvin} are non-local in time, and thus their numerical computation requires the storage of the time history of fields, which is very computationally expensive. To overcome this problem, the dispersion formula for lossless media, i.e., the  relation between the temporal frequency $\omega$ and spatial frequency $k$ ($\omega \approx c k$) with $c$ the sound speed, is used to replace the fractional time derivatives by fractional space
derivatives, which are non-local in space, rather than time \cite{Treeby-f}.
Using this method, the fractional time derivative is written as two fractional Laplacian operators in the form \cite{Treeby-a,Treeby-f}
\begin{align} \label{laplace}                                                                                                                                                                                                          \begin{split}                                                                                                                                                                                                                          \frac{\partial ^{y-1}}{\partial t^{y-1}} \approx & c^{y-1} \sin(\pi y/2) (- \nabla ^2 )^{(y-1)/2}\\
 &- c^{y-2} \cos (\pi y/2) (- \nabla ^2)^{(y-2)/2} \frac{\partial} {\partial t}.
\end{split}
\end{align}                                                                                                                                                                                                                            
This method provides a significant computational memory benefit, since at each time step the wavefield at all spatial positions is readily accessible. The system of viscoelastic wave equations enriched by fractional Laplacian operators can be used to describe absorption and physical dispersion behaviours over a wide range of frequencies and absorption values \cite{Treeby-f}. However, in elastic media since the compressional and shear waves travel at different speeds, separate dispersion relations must be considered for the compressional and shear parts of the wavefield. This requires that the particle velocity field is split into the compressional and shear components \cite{Treeby-f}. Throughout this work, superscripts $p$ and $s$ denote the compressional and shear parts of the fields, respectively. $v_i^p$ and $v_i^s$ are calculated in the form
\begin{align}  \label{spl}
\begin{split}
v_i^p & = q^p (v_i)  =  F^{-1} \Big\{  \hat{k_i}\hat{k_j}                F\big\{ v_j\big\}\Big\} \\
v_i^s & = q^s (v_i)  =  F^{-1} \Big\{ ( \delta_{ij}-\hat{k_i}\hat{k_j})  F\big\{ v_j\big\}\Big\},
\end{split}
\end{align}
where $F$ represents the Fourier transform operator, and $\hat{k}_i \hat{k}_j$ is the unit dyadic wavenumber tensor with $\hat{k}_i=k_i/k$ the unit vector in direction $i$ and $k=\left(\sum_i k_i^2 \right)^{1/2}$ the magnitude of wavenumber.
By splitting the particle velocity vector, the stress tensor is updated distinctly for compressional and shear parts in the form
\begin{align}  \label{stress:eq}                                                                                                                                                                                                       \begin{split}                                                                                                                                                                                                                          \frac{\partial \sigma_{ij}^{p,s}}{\partial t} &=  \lambda \left(\delta_{ij} \frac{\partial}{\partial x_l}  v_l^{p,s}\right) + \mu \left(\frac{\partial}{\partial x_j}v_i^{p,s}+\frac{\partial}{\partial x_i}v_j^{p,s}\right)+  \chi \left(\delta_{ij} \frac{\partial}{\partial x_l} \frac{\partial_{p,s}^{y-1}}{\partial t^{y-1}} v_l^{p,s} \right)\\
                                              &+  \eta \left(\frac{\partial}{\partial x_j}\frac{\partial_{p,s}^{y-1}}{\partial t^{y-1}}v_i^{p,s}+\frac{\partial}{\partial x_i}\frac{\partial_{p,s}^{y-1}}{\partial t^{y-1}}v_j^{p,s}\right).
\end{split}                                                                                                                                                                                                                            
\end{align}
Now, the fractional temporal derivatives in \eqref{stress:eq} can be replaced by fractional Laplacian operators using \eqref{laplace} with different sound speed maps for the compressional and shear waves. We will make this replacement in our forward model, and furthermore, following \cite{Treeby-f}, use the conservation of momentum to make the replacement
\begin{align} \label{velt:eq}
\frac{\partial v_i^{p,s}}{\partial t}=\frac{1}{\rho} \frac{\partial}{\partial x_j} \sigma_{ij}^{p,s}.
\end{align}
This is done to avoid having to compute time differences in the discretised model. To simplify the notation we introduce the operators
\begin{equation}
L_{c_{p,s}}^{y} = c_{p,s}^{y} (-\nabla^2)^{y/2}, \quad L_{c_{p,s}}^{y*} =  (-\nabla^2)^{y/2} c_{p,s}^{y},
\end{equation}
which are formal adjoints since $-\nabla^2$ is self-adjoint. Using this notation and the comments above, we will be using the following definition throughout the rest of this work including in \eqref{stress:eq}
\begin{equation} \label{laplace2}
\frac{\partial_{p,s}^{y-1}}{\partial t^{y-1}} v_i^{p,s} = \sin(\pi y/2) L_{c_{p,s}}^{y-1} v_i^{p,s} - \cos(\pi y/2) L_{c_{p,s}}^{y-2}  \frac{1}{\rho} \frac{\partial}{\partial x_j} \sigma_{ij}^{p,s}.
\end{equation}
Here, we assumed that $y$ is constant over the entire medium in the same way as \cite{Treeby-f}.

The continuous forward model for the wave propagation is completed with the conservation of momentum
\begin{align} \label{vel:eq}
\frac{\partial v_i}{\partial t}=\frac{1}{\rho}\sum_{p,s} \frac{\partial}{\partial x_j} \sigma_{ij}^{p,s}.
\end{align}

Equations \eqref{spl}, \eqref{stress:eq}, \eqref{laplace2}, and \eqref{vel:eq} together give a system of coupled partial differential equations which describe the propagation of PA waves in linear isotropic, lossy and heterogeneous viscoelastic media with an attenuation following the frequency power law.

\section{Continuous adjoint for viscoelastic wave equations} \label{methodi}
Let $\Omega \subset \mathbb{R}^d$ be a $d$-dimensional open, bounded set containing the initial pressure. We define the operator $\mathcal{G}$, which maps the compressional part of the stress tensor field to the pressure field in the form
\begin{align}      \label{fnc}
\mathcal{G}\sigma^p(x,t) = - \frac{1}{d} \delta_{ij} \sigma_{ij}^p(x,t) = p(x,t),
\end{align}
\noindent
that is the minus average trace of the compressional part of the stress tensor.
We also introduce $W(x,t) \in C_0^\infty (\Gamma \times \mathbb{R})$ for restricting the pressure $p(x,t)$ to the spatio-temporal field accessible to the sensors
with $\Gamma \subset \mathbb{R}^d$ an open, bounded set. Additionally, $\mathcal{M}$ maps the accessible part of the pressure field into the data measured by the sensors $\hat{P} \in \mathbb{R}^{N_s N_t}$ with $N_t,N_s \in \mathbb{N}$ the number of measurement time instants and the number of detectors, respectively.

\begin{definition} \label{def1}
Inspired by \cite{Arridge}, we define the PAT forward operator using the viscoelastic model in the form                                                                                                                                                                               
\begin{align}   \label{elast_map}
\begin{split}                                                                                                                                                                                                                          
&\Lambda: C_0^\infty(\Omega) \rightarrow \mathbb {R}^{N_s N_t}  \\
&\Lambda [p_0](x,t) =  \mathcal{M}  W(x,t) \mathcal{G} \sigma^p(x,t),
\end{split}                                                                                                                                                                                                                            
\end{align}
where $\sigma^{p,s}_{ij}$ and $v_i$ satisfy \eqref{stress:eq} and \eqref{vel:eq} with initial conditions
\begin{align}  \label{init1}
\sigma_{ij}^p(x,0) = -  \delta_{ij} p_0(x), \quad \sigma_{ij}^s(x,0) = 0,\quad v_i(x,0) = 0.
\end{align}
\end{definition}

\noindent
In the next  Lemma, we will calculate the adjoint of $\Lambda$ with respect to the $\mathcal{L}^2$ inner product, i.e.,
\begin{equation}
\Lambda^*:  \mathbb {R}^{N_s N_t} \rightarrow \mathcal{L}^2(\Omega).
\end{equation}
For this we also need the time reversal operator $\mathcal{R}$ defined by
\[
\mathcal{R}[p](x,t) = p(x,T-t).
\]

\begin{lemma}
The adjoint map $\Lambda^*$ can be calculated from $\Lambda^* [\hat{P}] (r,t) =  p_0^*(r)$, where $p_0^* = - \delta_{ij} \sigma_{ij}^{p*} (r,T)$, and $\sigma_{ij}^{p,s*}$ and $v_i^{{p,s}^*}$
satisfy the coupled equations
\begin{align} \label{vel2:eq}
\begin{split}
 & \rho \frac{{\partial v_i}^*}{\partial t}=  \sum_{p,s} q^{p,s}  \Bigg[ \Bigg( \frac{\partial}{\partial x_i} \Big (\lambda {\sigma_{ll}^{p,s}}^* \Big ) + 2 \frac{\partial}{\partial x_j} \Big( \mu {\sigma_{ij}^{p,s}}^* \Big)  \Bigg)\\
&  +  \sin(\pi y/2) L_{c_{p,s}}^{(y-1)*}  \Big)   \Bigg( \frac{\partial}{\partial x_i} \Big
(\chi {\sigma_{ll}^{p,s}}^* \Big ) + 2 \frac{\partial}{\partial x_j}  \Big( \eta {\sigma_{ij}^{p,s}}^*  \Bigg)  \Bigg]
\end{split}
\end{align}

\begin{align} \label{stress2:eq}
\begin{split}
{\textbf{v}_i^{p,s}}^*   &   =  v_i^*-\frac{1}{\rho} \cos(\pi y/2) L_{c_{p,s}}^{(y-2)*} \Bigg( \frac{\partial}{\partial x_i} \Big
(\chi {\sigma_{ll}^{p,s}}^* \Big ) + 2 \frac{\partial}{\partial x_j}  \Big( \eta {\sigma_{ij}^{p,s}}^* \Big)\Bigg) \\
\frac{\partial {\sigma_{ij}^p}^*}{\partial t} &= \frac{1}{2}\left (\frac{\partial {\textbf{v}^p_i}^*}{\partial x_j} + \frac{\partial {\textbf{v}^p_j}^*}{\partial x_i} \right )+(\mathcal{R} \mathcal{G} ^*W^*\mathcal{M}^* \hat{P}) \\
\frac{\partial {\sigma_{ij}^s}^*}{\partial t} &= \frac{1}{2}\left (\frac{\partial {\textbf{v}^s_i}^*}{\partial x_j} + \frac{\partial {\textbf{v}^s_j}^*}{\partial x_i} \right )
\end{split}
\end{align}
with initial conditions
\begin{align}  \label{init2}
 {\sigma_{ij}^{p,s}}^*(r,0) = 0, \quad v_i^*(r,0) = 0.
\end{align}

\begin{proof}
We will show that when $\sigma^{p,s}_{ij}$ and $v_i$ satisfy equations \eqref{stress:eq}, \eqref{vel:eq} and \eqref{init1}, and also $\sigma_{ij}^{{p,s}^*}$ and $v_i^*$ satisfy \eqref{vel2:eq}, \eqref{stress2:eq} and \eqref{init2}, then the forward map $\Lambda$ and adjoint $\Lambda^*$ must satisfy
\begin{align}    \label{adj_main}
\langle \Lambda [p_0] , \hat{P} \rangle_{\mathbb{R}^{N_s N_t}} =  \langle p_0  , \Lambda^*[\hat{P}]  \rangle_{\mathcal{L}^2(\mathbb{R}^d)}
\end{align}
for any $p_0 \in C_0^\infty (\Omega )$ and $\hat{P} \in \mathbb{R}^{N_s N_t}$.
Because it will make calculations easier, we first deal with the adjoint fields in a time reversed order (i.e. we make the change of variable $t \mapsto T-t$) and with $v_i^*$ replaced by $-v_i^*$ so that the initial conditions \eqref{init2} are actually final conditions $\sigma_{ij}^*(x,T)=0$ and $v_i^*(x,T)=0$. Accordingly, we have the following relation between $v_i$ and the adjoint field $v_i^*$
\[                                                                                                                                                                                                                                     \int_0^T \int_{\mathbb{R}^d} \rho \left ( \frac{\partial v_i}{\partial t} \overline{v_i^*} + v_i \overline{\frac{\partial v_i^*}{\partial t}} \right ) \ \mathrm{d} x \ \mathrm{d} t = 0.
\]                                                                                                                                                                                                                                     
Plugging \eqref{vel:eq} into the first integrand in the above equation gives                                                                                                                                                                                          
\begin{equation} \label{ze}                                                                                                                                                                                                                                   \int_0^T \int_{\mathbb{R}^d} \sum_{p,s}  \frac{\partial \sigma_{ij}^{p,s}}{\partial x_j} \overline{v_i^*} + \rho v_i \overline{\frac{\partial v_i^*}{\partial t}}  \ \mathrm{d} x \ \mathrm{d} t = 0.
\end{equation}                                                                                                                                                                                                                                     
Integrating-by-parts, and using that $\sigma^{p,s}_{ij} \rightarrow 0$ at infinity, we end up with
\begin{equation}\label{1}                                                                                                                                                                                                              \int_0^T \int_{\mathbb{R}^d} \sum_{p,s} - \sigma_{ij}^{p,s}\frac{1}{2} \overline{\left ( \frac{\partial v_i^*}{\partial x_j} + \frac{\partial v_j^*}{\partial x_i}\right )} + \rho v_i \overline{\frac{\partial v_i^*}{\partial t}}  \ \mathrm{d} x \ \mathrm{d} t = 0.
\end{equation}                                                                                                                                                                                                                         
In the above equation, we also used the symmetry of the stress tensor $\sigma_{ij}^{p,s} = \sigma_{ji}^{p,s}$. Now we apply the same procedure to the stress tensor. Using the final conditions $\sigma_{ij}^{p,s*}(x,T) = 0$ yields
\begin{align}  \label{a1}                                                                                                                                                                                                                                   \int_0^T \int_{\mathbb{R}^d}    \frac{\partial \sigma_{ij}^p}{\partial t} \overline{{\sigma_{ij}^p}^*} + \sigma_{ij}^p \overline{\frac{\partial {\sigma_{ij}^p}^*}{\partial t}}  \ \mathrm{d} x \ \mathrm{d} t =  \int_{\mathbb{R}^3} p_0\ \overline{{\sigma_{ii}^p}^*(x,0)} \ \mathrm{d} x  = -\langle p_0, \Lambda^*[\hat{P}] \rangle_{\mathcal{L}^2(\mathbb{R}^d)}.                                                                                                                                                                                                                                   \end{align}                                                                                                                                                                                                                                     
and
\begin{align}    \label{a2}                                                                                                                                                                                                                                \int_0^T \int_{\mathbb{R}^d}    \frac{\partial \sigma_{ij}^s}{\partial t} \overline{{\sigma_{ij}^s}^*} + \sigma_{ij}^s \overline{\frac{\partial {\sigma_{ij}^s}^*}{\partial t}}  \ \mathrm{d} x \ \mathrm{d} t =  0.                                                                                                                                                                                                                                     \end{align}
Now, plugging \eqref{stress:eq} into the first integrands in the left-hand sides of \eqref{a1} and \eqref{a2} and then adding these two equations results in
\[                                                                                                                                                                                                                                     \begin{split}                                                                                                                                                                                                                          \int_0^T \int_{\mathbb{R}^d} & \sum_{p,s} \Bigg[ \bigg ( \lambda \delta_{ij} \frac{\partial v_l^{p,s}}{\partial x_l} + \mu \Big ( \frac{\partial v_i^{p,s}}{\partial x_j} + \frac{\partial v_j^{p,s}}{\partial x_i}\Big) +\chi \Big( \delta_{ij}
\frac{\partial}{\partial x_l} \frac{\partial_{p,s}^{y-1}}{\partial t^{y-1}}  v_l^{p,s}  \Big)\\
& + \eta \Big( \frac{\partial}{\partial x_j} \frac{\partial_{p,s}^{y-1}}{\partial t^{y-1}} v_i^{p,s}+\frac{\partial}{\partial x_i} \frac{\partial_{p,s}^{y-1}}{\partial
t^{y-1}}v_j^{p,s}\Big)\bigg )                                                                                                                                                                                                             \overline{{\sigma_{ij}^{p,s}}^*}                                                                                                                                                                                                             + \sigma_{ij}^{p,s} \overline{\frac{\partial {\sigma_{ij}^{p,s}}^*}{\partial t}} \Bigg] \ \mathrm{d} x \ \mathrm{d} t                                                                                                                                             =   -\langle p_0, \Lambda^*[\hat{P}] \rangle_{\mathcal{L}^2(\mathbb{R}^d)}.                                                                                                                                                                         \end{split}                                                                                                                                                                                                                            \]                                                                                                                                                                                                                                     
Taking integration-by-parts to the first term in the bracket in the above equation, together with the fact that ${\sigma_{ij}}^*$ is symmetric, gives                                                                                                                                                                                                      
\begin{equation}\label{2}                                                                                                                                                                                                              \begin{split}                                                                                                                                                                                                                          \int_0^T \int_{\mathbb{R}^d}  & \sum_{p,s} \Bigg[ -\Bigg(\overline{\frac{\partial}{\partial x_i} \Big (\lambda {\sigma_{ll}^{p,s}}^* \Big ) + 2 \frac{\partial}{\partial x_j} \Big ( \mu {\sigma_{ij}^{p,s}}^* \Big ) + \frac{\partial}{\partial x_i} \Big
(\chi {\sigma_{ll}^{p,s}}^* \Big )}\frac{\partial_{p,s}^{y-1}}{\partial t^{y-1}} \\
&+ 2 \overline{\frac{\partial}{\partial x_j}  \Big( \eta {\sigma_{ij}^{p,s}}^*\Big)}\frac{\partial_{p,s} ^{y-1}}{\partial t^{y-1}} \Bigg) v_i^{p,s}
+ \sigma_{ij}^{p,s} \overline{\frac{\partial {\sigma_{ij}^{p,s}}^*}{\partial t}} \Bigg] \ \mathrm{d} x \ \mathrm{d} t =-\langle p_0, \Lambda^*[\hat{P}] \rangle_{\mathcal{L}^2(\mathbb{R}^d)}.                                                                                                                                                                           \end{split}                                                                                                                                                                                                                            \end{equation}                                                                                                                                                                                                                         
Now, adding equations \eqref{1} and \eqref{2} yields                                                                                                                                                                                   
\begin{equation} \label{22}
 \begin{split}                                                                                                                                                                                                                          \int_0^T \int_{\mathbb{R}^d}  &  \rho v_i \overline{ \frac{\partial v_i^*}{\partial t}} - \sum_{p,s} \Big[ \overline{\frac{\partial}{\partial x_i} \Big (\lambda {\sigma_{ll}^{p,s}}^* \Big ) + 2\frac{\partial}{\partial x_j} \Big ( \mu {\sigma_{ij}^{p,s}}^*\Big )} \\
& +\overline{ \Bigg (\frac{\partial}{\partial x_i} \Big (\chi {\sigma_{ll}^{p,s}}^* \Big ) + 2 \frac{\partial}{\partial x_j}  \Big( \eta {\sigma_{ij}^{p,s}}^* \Big ) \Bigg )} \frac{\partial_{p,s}^{y-1}}{\partial t^{y-1}}  \Big] v_i^{p,s}    \\
+& \sum_{p,s} \overline{ \bigg[   \frac{\partial {\sigma_{ij}^{p,s}}^*}{\partial t} - \frac{1}{2}\left (\frac{\partial v_i^*}{\partial x_j} + \frac{\partial v_j^*}{\partial x_i} \right ) \bigg]} \sigma_{ij}^{p,s}    \ \mathrm{d} x \ \mathrm{d} t  = -\langle p_0, \Lambda^*[\hat{P}] \rangle_{\mathcal{L}^2(\mathbb{R}^d)}.                                                                                                                                                                           \end{split}                                                                                                                                                                                                                            \end{equation}
Now, plugging the fractional Laplacian operators defined in equation \eqref{laplace2} into the second line in equation \eqref{22}, together with \eqref{spl} and using \eqref{velt:eq} yields
\begin{align}  \label{3}
\begin{split}
 \int_0^T & \int_{\mathbb{R}^d}   \Bigg( \rho \overline{\frac{\partial v_i^*}{\partial t}} -  \sum_{p,s}  \Bigg[ \overline{\frac{\partial}{\partial x_i} \Big (\lambda {\sigma_{ll}^{p,s}}^* \Big ) + 2 \frac{\partial}{\partial x_j} \Big( \mu {\sigma_{ij}^{p,s}}^* \Big)} \\
+ & \overline{\Bigg( \frac{\partial}{\partial x_i} \Big
(\chi {\sigma_{ll}^{p,s}}^* \Big ) + 2 \frac{\partial}{\partial x_j}  \Big( \eta {\sigma_{ij}^{p,s}}^* \Big)\Bigg)}  \sin(\pi y/2) L_{c_{p,s}}^{y-1}     \Bigg] q^{p,s} \Bigg) v_i \\
+ & \sum_{p,s} \overline{\Bigg( \frac{\partial}{\partial x_i} \Big
(\chi {\sigma_{ll}^{p,s}}^* \Big ) + 2 \frac{\partial}{\partial x_j}  \Big( \eta {\sigma_{ij}^{p,s}}^* \Big)\Bigg)}  \cos(\pi y/2) L_{c_{p,s}}^{y-2}   \frac{1}{\rho}\frac{\partial}{\partial x_j} \sigma_{ij}^{p,s}  \\
+& \sum_{p,s} \overline{ \bigg[   \frac{\partial {\sigma_{ij}^{p,s}}^*}{\partial t} - \frac{1}{2}\left (\frac{\partial v_i^*}{\partial x_j} + \frac{\partial v_j
^*}{\partial x_i} \right ) \bigg]} \sigma_{ij}^{p,s}
\ \mathrm{d} x \ \mathrm{d} t= -\langle p_0, \Lambda^*[\hat{P}] \rangle_{\mathcal{L}^2(\mathbb{R}^d)}.
\end{split}
\end{align}
In the above equation, we also used the linearity of operator $q^{p,s}$ with respect to $v_i$. By taking integration-by-parts to the third line and using the symmetry of $\sigma_{ij}^{p,s}$ in the same way as the first integrand in \eqref{ze}, we can see from \eqref{stress2:eq} that if the integral in the first two lines of \eqref{3} is equal to zero, then \eqref{adj_main} holds and the proof is complete (recall again that relative to \eqref{vel2:eq}, we have reversed the time, and changed $v_i^*$ to $-v_i^*$). So we now focus on the first two lines of \eqref{3} which we will denote $\mathcal{I}$. Considering that $q^{p,s}$ is self-adjoint yields

\begin{align}  \label{5}
\begin{split}
\mathcal{I} = & \int_0^T  \int_{\mathbb{R}^d}  \overline{ \Bigg[ \rho \frac{\partial v_i^*}{\partial t} -  \sum_{p,s} q^{p,s} \bigg[ \bigg( \frac{\partial}{\partial x_i} \Big (\lambda {\sigma_{ll}^{p,s}}^* \Big ) + 2 \frac{\partial}{\partial x_j} \Big( \mu {\sigma_{ij}^{p,s}}^* \Big) \bigg)}   \\
 & \overline{+ \sin(\pi y/2) L_{c_{p,s}}^{(y-1)*}   \Bigg( \frac{\partial}{\partial x_i} \Big
(\chi {\sigma_{ll}^{p,s}}^* \Big ) + 2 \frac{\partial_{p,s}}{\partial x_j}  \Big( \eta {\sigma_{ij}^{p,s}}^* \Big)\Bigg) \bigg]  \Bigg]} v_i   \ \mathrm{d} x \ \mathrm{d} t.
\end{split}
\end{align}
Since $v_{i}^*$ and $\sigma_{ij}^{p,s*}$ satisfy \eqref{vel2:eq} we can now see that in fact $\mathcal{I} = 0$, and so the proof is complete.
%
%
\end{proof}
\end{lemma}
\noindent
Setting $\chi,\eta=0$ and ignoring the splitting operator results in the general form of the adjoint for lossless media, which can be adapted to other attenuation models (see for example \cite{Kowar,Kowar-b}).
\noindent

\section{Numerical computation} \label{nmethodf}

Having found an analytically exact method of computing the adjoint operator $\Lambda^*$ in the previous section, we now consider in more detail the discretisation and computation of the forward operator $\Lambda$ and adjoint operator $\Lambda^*$. We approximate the fields on a uniform rectilinear grid staggered in space and time \cite{Tabei,Firouzi}. We denote the position of a given grid point in Cartesian coordinates by $x_\zeta$ where $\zeta = (\zeta_1, \ ... \ , \zeta_d) \in\{1, \ ... \ , N_1\} \times \ ... \ \times \{1, \ ... \ , N_d\}$ with $N=\prod_{i=1}^d N_i$ the total number of grid points along $d$ dimensions. The grid spacing along the $i$th direction will be denoted $\Delta x_{i}$. Let also $n$ denote the iteration corresponding to time $t_n= n\Delta t$ with $n \in \{-1, ...,N_t-1\}$. Using a staggered temporal grid, $t_{n+1/2}= n \Delta t+\Delta t/2 $.

We discretise the spatial derivatives by a pseudo-spectral method. The k-space correction is also applied to the spatial derivatives in order to minimize the numerical dispersion errors due to the time integration. Using a staggered spatial grid, these give the spatial gradient in direction $i$ in the form
\begin{align}  \label{grf}
\frac{\partial_{p,s}\left\{ \cdot \right\} }{\partial x_i^{\pm}}=F^{-1}\left\{\textbf{i}k_i \sinc(c_0^{p,s} k\Delta t/2)e^{\pm\textbf{i} k_i \Delta x_i/2}F\left\{\cdot\right\} \right\},
\end{align}
where, as opposed to \eqref{spl}, $F$ and $F^{-1}$ denote the discrete Fourier transform and its inverse, while $c_0^{p,s}$ is the reference sound speed associated with the compressional and shear parts of the fields. The reader is referred to \cite{Tabei,Firouzi} for further details on the k-space pseudo-spectral method.

To avoid spurious reflections at the boundaries, it is necessary that the outward travelling waves that reach the edge of the domain are absorbed by perfectly matched layers (PMLs) \cite{Firouzi,Mitsuhashi}. Using PMLs, the general evolution equation $\frac{\partial R(x,t)}{\partial t}=\beta(x,t)$ is transformed into the form \cite{Tabei}
\begin{align} \label{PML}
\frac{\partial R(x,t)}{\partial t}+\alpha_a R(x,t)=\beta(x,t),
\end{align}
where $\alpha_a$ is the attenuation coefficient associated with the PML, which is tapered within the PML thickness at each side of the grid (cf. \cite{Tabei}, Eq.(27)). This yields
\begin{align}
\frac{\partial \left(e^{\alpha_a t}R(x,t)\right)}{ \partial t}= e^{\alpha_a t} \beta (x,t).
\end{align}
Using a staggered temporal grid, this is approximated as
\begin{align}  \label{PML-st}
\frac{e^{\alpha_a(t+\Delta t)}R(x_\zeta,t+\Delta t)- e^{\alpha_a t} R(x_\zeta,t)}{\Delta t}=e^{\alpha_a (t+\Delta t/2)}\beta(x_\zeta,t+\Delta t/2).
\end{align}
This gives the update
\begin{align}
R(x_\zeta,t+\Delta t)=e^{-\alpha_a \Delta t/2} \Big[ e^{-\alpha_a \Delta t/2} R(x_\zeta,t) + \Delta t \beta(x_\zeta,t+\Delta t/2) \Big].
\end{align}
Using direction-dependent PMLs, the field variables are split into directions along the Cartesian coordinates $m \in \{1,...,d\}$ \cite{Firouzi}. In the sequel, the directions associated with PMLs are written to the left of the fields. We define the diagonal PML attenuation matrices $A_m \in \mathbb{R}^{N \times N}$ by
\begin{align}  \label{A}
A_m= \mathrm{diag}(e^{- _{m}{\alpha}_a \Delta t/2}).
\end{align}
Note that $-_{m}\alpha_a$ depends on the grid point here.

To accommodate the staggered grid we introduce the operators $\mathbb{T}_{i}$ which shift the point $x$ by $\Delta x_i/2$ in the $i$th coordinate, i.e., $x_i$ changes to $x_i + \Delta x_i/2$. We will also use the same notation for the corresponding operator acting on functions defined by
\begin{align}
\mathbb{T}_i f(x) = f (\mathbb{T}_i x).
\end{align}
The discretised particle velocity vector field is denoted by $_{m}v_{(i;\zeta;n)} \in \mathbb{R}^d_m \times \mathbb{R}^d_i \times \mathbb{R}_\zeta^N \times \mathbb{R}^{N_t+1}_n$ and is approximated on a staggered spatial grid as
\begin{align}
_{m}v_{(i;\zeta;n)} \approx \ \mathbb{T}_{i}\ {_{m}v_{i}} \left(  x_\zeta ,t_n  \right) .
\end{align}
The $p$ and $s$ parts of the discretised stress tensor field are denoted by $_{m}\sigma_{(ij;\zeta;n)}^{p,s} \in \mathbb{R}^d_m \times \mathbb{R}^d_i \times\mathbb{R}^d_j \times \mathbb{R}_\zeta^N \times \mathbb{R}^{N_t+1}_n$ and are approximated on a staggered grid as
\begin{align}
_{m}\sigma_{(ij;\zeta;n)}^{p,s} \approx
\begin{cases}
{_{m}\sigma}_{ij}^{p,s}\left(  x_\zeta,t_n  \right) & \text{if} \quad i = j\\
\mathbb{T}_i \mathbb{T}_j\ {_{m}\sigma}_{ij}^{p,s}\left(  x_\zeta,t_n  \right)  & \text{if} \quad i \neq j.
\end{cases}
\end{align}

\noindent
Because of using a staggered grid, the unit dyadic tensor in \eqref{spl} will be in the form \cite{Firouzi}
\begin{align} \label{si1}
(\hat{k}_{i}\hat{k}_{j})_{\text{staggered}} = (\hat{k}_{i}\hat{k}_{j})_{\text{nonstaggered}} \times \xi_{ij},
\end{align}
where
\begin{align}  \label{si2}
\xi_{ij}= e^{+\textbf{i}\left( k_ {i}\Delta x_{i}-k_{j}\Delta x_{j} \right)/2}                                                                                                                                                                                                                         
\end{align}
is the shifting operator with $\textbf{i}$ standing for the imaginary number \cite{Firouzi}.

\noindent
We also define
\begin{align} \label{lp}
\tau_\text{dis}^{p,s} & = C_{p,s}^{y-1} \sin{(\pi y/2)}    \\
\tau_\text{abs}^{p,s} & = C_{p,s}^{y-2} \cos{(\pi y/2)}
\end{align}
with $C_{p,s} \in \mathbb{R}^N$ the discretised form of $c_{p,s}$. Using a staggered grid, we define the medium's parameters as diagonal matrices of size $N \times N$ in the form
\begin{align} \label{mp}
\begin{split}
\bar{\rho}_i                       & =\diag \left(\mathbb{T}_i  \rho\right)\\
\bar{\lambda}                & =\diag \left( \lambda \right)\\
\bar{\mu}_{ij}                & = \begin{cases}
\diag \left(\mu \right) & \text{if} \quad i = j\\
\diag \left( \mathbb{T}_i \mathbb{T}_j\mu \right )  & \text{if} \quad i \neq j.
\end{cases}\\
\bar{\chi}                       & =\diag \left( \chi \right)\\
\bar{\eta}_{ij}               & = \begin{cases}
\diag \left( \eta \right) & \text{if} \quad i = j\\
\diag \left( \mathbb{T}_i \mathbb{T}_j\eta \right )  & \text{if} \quad i \neq j.
\end{cases}\\
\bar{\tau}_{i,\text{dis}}     &= \diag  \left(\mathbb{T}_i \tau_\text{dis} \right)\\
\bar{\tau}_{i,\text{abs}}   &= \diag  \left(\mathbb{T}_i \tau_\text{abs} \right)
\end{split}
\end{align}
where $\rho, \lambda, \mu,\chi$, and $\eta$ on the right hand sides in \eqref{mp} are the medium parameters evaluated at the $N$ grid points. In the formulas that follow for the discretised model, these matrices are always understood to act on discretised fields in the index $\zeta$ corresponding to the spatial grid.

We also introduce the $N \times N$ matrices discretising the relevant fractional Laplacian operators as
\begin{align} \label{fracLap}
\begin{split}
\bar{Y}_{\text{dis}} &=   F^{-1}  \Big \{  k^{y-1}  F\{\cdot \}\Big \} \\
\bar{Y}_{\text{abs}} &=   F^{-1}  \Big \{  k^{y-2}  F\{\cdot \}\Big \}.
\end{split}
\end{align}
Finally, we define the following function which we will use to simplify some of the formulas
\begin{align}
\begin{split}
h(i,j)=
\begin{cases}
 +1  &\text{if}  \quad   i=j  \\
 -1   &\text{if}  \quad   i\neq j.
\end{cases}
\end{split}
\end{align}


\subsection{Forward model.}   \label{f}
In the sequel, the approximation of the system of viscoelastic wave equations defined by equations \eqref{stress:eq}, \eqref{vel:eq} and \eqref{init1} based on the details given above will be outlined. A code is available in \textit{k-Wave} toolbox for describing wave propagation in heterogeneous elastic media using the k-space pseudo-spectral method \cite{treeby-b,Treeby-tool}. We enriched this code by the splitting operator and fractional Laplacians in order to include absorption and physical dispersion following the frequency power law. This code is outlined as follows.

While in the continuous model we assume the initial pressure is instantaneous, in the discretised model we introduce the initial pressure at $t=0$ to the forward model as an additive source split over the time interval $t\in \big[-\Delta t/2,+\Delta t/2\big]$. For this, $p(t=0)=1$ is approximated as $P(n=\big[-1/2 \ +1/2\big])=\frac{1}{\Delta t}\big[0.5 \ 0.5\big]$ (cf. \cite{Arridge}, Appendix B). Considering this, together with \eqref{init1} and dividing the source by PML directions, gives a source in the form
\begin{align}    \label{s1}
_{m}s_{(ij;\zeta;n+1/2)} =
\begin{cases}
- \frac{\delta_{ij}}{2d \Delta t} \mathbb{S}P_0  \quad   &n=-1,0 \\
            0                                              \quad   &\text{otherwise},
\end{cases}
\end{align}
where $P_0$ denotes the discretised form of $p_0(r,t)$, and $\mathbb{S}$ is a symmetric smoothing operator that is used for mitigating unexpected oscillations in propagation of the initial pressure $P_0$ (For further details, the reader is referred to \cite{Arridge}, Appendix B).
\smallskip

\noindent \textit{Start at iterate $n=-1$ with initial conditions $_{m}\sigma^{p,s}_{(ij;\zeta;n=-1)}=0$ and $_{m}v_{(i;\zeta,n=-3/2)}=0$, and terminate at iterate $n=N_t-2$.}
\smallskip

\noindent
\textit{1. Update the particle velocity field:}
\begin{align} \label{num-vel}
\begin{split}
_{m}v_{(i;\zeta;n+\frac{1}{2})}&= A_m   \bigg[A_m \  {_{m}v_{(i;\zeta;n-\frac{1}{2})}}
+ \frac{\Delta t}{\bar{\rho}_i}\sum_{p,s} \frac{ \partial_{p,s}}{\partial x_m^{h(i,m)} }    \sigma_{(im;\zeta;n)}^{p,s}   \bigg]  \\
v_{(i;\zeta;n+\frac{1}{2})}&=\sum_{m=1}^d {_{m}v_{(i;\zeta;n+\frac{1}{2})}}.
\end{split}
\end{align}

\noindent
\textit{2. Split the particle velocity field into compressional and shear parts:}
\begin{align}   \label{spl2}
v^{p,s}_{(i;\zeta;n+\frac{1}{2})}= Q^{p,s} v_{(i;\zeta;n+\frac{1}{2})}.
\end{align}
Here, $ Q^{p,s}$ denotes the discretised form of functions $q^{p,s}$ defined in \eqref{spl}.

\noindent
\textit{3. Update the stress tensor field:}
\begin{align} \label{num-s}
\begin{split}
_{m}\sigma_{(ij;\zeta;n+1)}^{p,s}&= A_m   \Bigg[A_m \  {_{m}\sigma_{(ij;\zeta;n)}^{p,s}}+\Delta t    \  \bigg[  \bar{\lambda}  \delta_{ij} \frac{ \partial_{p,s}}{\partial x_m^- }   v_{(m;\zeta;n+\frac{1}{2})}^{p,s}   \\
&+ \bar{\mu}_{ij} \bigg( \delta_{mj} \frac{ \partial_{p,s}}{\partial x_j^{-h(i,j)} }   v_{(i;\zeta;n+\frac{1}{2})}^{p,s}+
 \delta_{mi}  \frac{\partial_{p,s}}{\partial x_i^{-h(i,j)} }  v_{(j;\zeta;n+\frac{1}{2})}^{p,s} \bigg) \\
&+   \bar{\chi}  \delta_{ij} \frac{  \partial_{p,s}}{\partial x_m^- }  \left(\bar{\tau}_{m,\text{dis}}^{p,s}    \  \bar{Y}_{\text{dis}}\ v_{(m;\zeta;n+\frac{1}{2})}^{p,s} - \bar{\tau}_{m,\text{abs}}^{p,s}\  \bar{Y}_{\text{abs}}\     \partial_t v_{(m;\zeta;n+1/2)}^{p,s} \right)    \\
&+  \bar{\eta}_{ij} \delta_{mj} \frac{ \partial_{p,s}}{\partial x_j^{-h(i,j)} }  \left(\bar{\tau}_{i,\text{dis}}^{p,s} \  \bar{Y}_{\text{dis}}\  v_{(i;\zeta;n+\frac{1}{2})}^{p,s} -  \bar{\tau}_{i,\text{abs}}^{p,s}\  \bar{Y}_{\text{abs}}\ \partial_t v_{(i;\zeta;n+1/2)}^{p,s}\right) \\
&+ \bar{\eta}_{ij} \delta_{mi} \frac{ \partial_{p,s}}{\partial x_i^{-h(i,j)} } \left(\bar{\tau}_{j,\text{dis}}^{p,s}    \  \bar{Y}_{\text{dis}}\  v_{(j;\zeta;n+\frac{1}{2})}^{p,s} - \bar{\tau}_{j,\text{abs}}^{p,s}\  \bar{Y}_{\text{abs}}\ \partial_t v_{(j;\zeta;n+1/2)}^{p,s}\right)  \bigg]\Bigg].
\end{split}
\end{align}
where
\begin{align}  \label{num-s2}
\partial_t v_{(i;\zeta;n+1/2)}^{p,s} = \sum_{m=1}^d A_m \frac{1}{\bar{\rho}_i} \frac{\partial_{p,s}}{\partial x_m^{h(i,m)} }    \sigma_{(im;\zeta;n)}^{p,s}
\end{align}

\noindent
\textit{4. Add source:}
\begin{align}      \label{s2}
\begin{split}
_{m}\sigma_{(ij;\zeta;n+1)}^{p}& \leftarrow{_{m}\sigma_{(ij;\zeta;n+1)}^{p}} + \Delta t \ {_{m}s_{(ij;\zeta;n+1/2)}} \\
\sigma_{(ij;\zeta;n+1)}^{p,s}&=\sum_{m=1}^d {_{m}\sigma_{(ij;\zeta;n+1)}^{p,s}}.
\end{split}
\end{align}

\noindent
\textit{5. Compute the pressure field and map it to detected data at ultrasound detectors:}\\
We use $\mathbb{G}$,$\mathbb{W}$ and $\mathbb{M}$ for denoting the discretised variants of $\mathcal{G}$,$W$ and $\mathcal{M}$. Correspondingly, at each iterate the pressure field is computed by
\begin{align}    \label{G}
p_{(\zeta;n+1)}= \mathbb{G} \sigma^p_{(ij;\zeta;n+1)}= -\frac{1}{d} \sum_{l,m=1}^d {_{m}\sigma_{(ll;\zeta;n+1)}^p},
\end{align}
and is then interpolated to ultrasound detectors using trilinear interpolation \cite{Arridge,Mitsuhashi}, i.e.,
\begin{align}\label{M}
\hat{P}_{n+1}= \mathbb{M} \mathbb{W} p_{(\zeta;n+1)}  ,
\end{align}
where $\mathbb{M} \mathbb{W} \in \mathbb{R}^{N_s \times N}$ is a map from the pressure at grid points to the pressure at the detector positions, and $\hat{P}_n \in \mathbb{R}^{N_s}$ is the vector of measured pressure data at iteration $n=0,...,N_t-1$.

\subsection{Analytic adjoint model} \label{nmethodi}

\noindent
The continuous adjoint model defined by equations \eqref{vel2:eq}, \eqref{stress2:eq} and \eqref{init2} are solved numerically as follows. For brevity, we ignore the superscript $*$ for denoting the adjoint fields in the discretised case. Before defining the time stepping procedure for the adjoint, we first define the additive source. To account for splitting of $P_0$ over the first two temporal iterations in the forward model (cf. equation \eqref{s1}), we define the order reversed adjoint measured data in the form \cite{Arridge}
\begin{align}      \label{adj_data}
\begin{split}
\hat{P}_{n+1/2}^\text{adj}= \frac{1} {2 \Delta t}
\begin{cases}
\hat{P}_{N_t-1},                                                   \quad &n=-1  \\
\hat{P}_{N_t-n-1}+\hat{P}_ {N_t-n-2},                              \quad &n=0,...,N_t-2 \\
\hat{P}_0,                                                        \quad &n=N_t-1
\end{cases}.
\end{split}
\end{align}
The adjoint measured data is mapped from ultrasound detector positions to an additive source that is defined at grid points using
\begin{align}
_{m}s_{(ij;\zeta;n+1/2)}=\mathbb{G}^T \mathbb{W}^T\mathbb{M}^T \hat{P}_{n+1/2}^\text{adj},
\end{align}
where
\begin{align}
\mathbb{G}^T  =   -  1_m \otimes  \frac{\delta_{ij}}{d}.
\end{align}

\noindent
\textit{Start at iterate $n=-1$ with initial conditions $_{m}\sigma_{(ij;\zeta;n=-1)}=0$ and $_{m}v_{(i;\zeta,n=-3/2)}=0$, and terminate at iterate $n=N_t-2$.}

%

\noindent
\textit{1. Update the particle velocity field:}
\begin{align} \label{cvx}
\begin{split}
_{m}v_{(i;\zeta;n+\frac{1}{2})}^{p,s}&= A_m   \Bigg[A_m \  {_{m}v_{(i;\zeta;n-\frac{1}{2})}}\\
&+ \frac{\Delta t}{\bar{\rho}_i}  \bigg[ \sum_{p,s} \sum_{j=1}^d  Q^{p,s}\Big[ \frac{\partial_{p,s}}{\partial x_i^{+} }  \big (\bar{\lambda} \  {_{i}\sigma_{(jj;\zeta;n)}^{p,s}} \big) +2 \frac{ \partial_{p,s}}{\partial x_j^{h(i,j)} } ( \bar{\mu}_{ij} \ {_{j}\sigma_{(ij;\zeta;n)}^{p,s}} ) \\
&+ \bar{Y}_{\text{dis}}\  \bar{\tau}_{i,\text{dis}}^{p,s} \Bigg(
 \frac{ \partial_{p,s}}{\partial x_i^{+} }  \big( \bar{\chi} \ {_{i}\sigma_{(jj;\zeta;n)}^{p,s}}\big )+ 2 \frac{ \partial_{p,s}}{\partial x_j^{h(i,j)} } (  \bar{\eta}_{ij} \ {_{j}\sigma_{(ij;\zeta;n)}^{p,s}}) \Bigg)\Big]
 \bigg] \Bigg]
\end{split}
\end{align}

\textit{2. Add the absorption term to the particle velocity field:}
\begin{align} \label{cvx2}
\begin{split}
_{m}\textbf{v}_{(i;\zeta;n+\frac{1}{2})}^{p,s}&= _{m}v_{(i;\zeta;n+\frac{1}{2})}\\
&-A_m \frac{1}{\bar{\rho}_i}\sum_{j=1}^d  \bar{Y}_{\text{abs}}\  \bar{\tau}_{i,\text{abs}}^{p,s}  \Bigg(
\frac{ \partial_{p,s}}{\partial x_i^{+} } \big (  \bar{\chi} \  {_{i}\sigma_{(jj;\zeta;n)}^{p,s}}\big ) + 2 \frac{ \partial_{p,s}}{\partial x_j^{h(i,j)} } (   \bar{\eta}_{ij} \ {_{j}\sigma_{(ij;\zeta;n)}^{p,s}}) \Bigg)
\end{split}
\end{align}

\noindent
Note that we are not using the summation convention in these formulas.

\noindent
\textit{3. Update the stress tensor field:}
\begin{align}   \label{cs:xx}
\begin{split}
_{m}\sigma_{(ij;\zeta;n+1)}= A_m  & \bigg[A_m \  {_{m}\sigma_{(ij;\zeta;n)}}+ \frac{\Delta t}{2} \Big(   \frac{ \partial_{p,s}}{\partial x_j^{-h(i,j)} }    {_{j}\textbf{v}_{(i;\zeta;n+\frac{1}{2})}^{p,s}}
+          \frac{\partial_{p,s}}{\partial x_i^{-h(i,j)} }    {_{i}\textbf{v}_{(j;\zeta;n+\frac{1}{2})}^{p,s}}
\Big)  \bigg]
\end{split}
\end{align}

\noindent
\textit{4. Add source:}

\begin{align}    \label{source_adj}
\begin{split}
_{m}\sigma_{(ij;\zeta;n+1)}^{p}&\leftarrow {_{m}\sigma_{(ij;\zeta;n+1)}^{p}}+    \Delta t \  {_{m}s_{(ij;\zeta;n+1/2)}}
\end{split}
\end{align}

\noindent
\textit{5. Compute the pressure field at final iterate and apply smoothing}:\\
\begin{align}
p_{(\zeta;n=N_t-1)}=-\mathbb{S} \left(  \frac{1}{d}  \sum_{l,m=1}^d  {_{m}\sigma_{(ll;\zeta;n=N_t-1)}^p} \right).
\end{align}

\section{Adjoint for discretised viscoelastic forward model}   \label{s6}

In this section, we will calculate the adjoint of the viscoelastic forward model $\Lambda$ defined by \eqref{stress:eq}, \eqref{vel:eq} and initial conditions in \eqref{init1} based on the \textit{discretize-then-adjoint} method. To do this, we consider the discretised equations \eqref{num-vel} and \eqref{num-s} in a matrix form. Accordingly, let the particle velocity vector $\bar{v}_{n-1/2} \in \mathbb{R}^{Nd^2}$ at each time step be made up of the components $_{m}v_{(i;\zeta;n-1/2)} \in \mathbb{R}^N \ (i,m \in \{1,...,d\})$. Let us also define the stress tensor as vector $\bar{\sigma}_n \in  \mathbb{R}^{12N} (d=2)$ or $\bar{\sigma}_n \in \mathbb{R}^{30N} (d=3)$ composed of the components $_{m}\sigma_{(ij;\zeta,n)}^{p,s} \in \mathbb{R}^N \ ( i,j,m \in \{1,...,d\})$. Note that for the latter, we used the symmetry of $\sigma_{ij}$, together with the fact that $_{m}\sigma_{ij}=0$ if $m\notin \{i,j\}$ \cite{Firouzi} to reduce the number of degrees of freedom.
We also define $X_n \in R^{39N}$ (3D case) as a stack of the particle velocity and stress fields at times corresponding to time step $n$ in the form $X_n=\big[\left(\bar{v}_{n-1/2}\right)^T \ \left(\bar{\sigma}_n\right)^T\big]^T$. Let also $T_0 \in \mathbb{R}^{39NN_t \times N}$ give the map from the discretised initial pressure $P_0 \in \mathbb{R}^N$ to an additive source (cf. \eqref{s1}), which we will write as
\begin{align} \label{dip}
S =\mathcal{S} P_0 \in \mathbb{R}^{39NN_t}.
\end{align}
We will also write $S_{n+1/2} = \mathcal{S}_n P_0$ for the source at time step $n$. In particular $\mathcal{S}_n = 0$ except when $n = -1$ or $0$.
The time sequence of fields at steps $ (n \in \{-1,..., N_t-2\})$ is then given by
\begin{align}   \label{t}
X_{n+1}=T X_{n} + S_{n+1/2},
\end{align}
where $T \in \mathbb{R}^{39N \times 39N}$ implements \eqref{num-vel} and \eqref{num-s}, and $X_{-1} = 0$ (cf. section \ref{f}). Here for brevity the operators $\mathcal{S}_n$ are given using \eqref{s1} multiplied by $\Delta t$, and thus multiplication by $\Delta t$ is neglected in the second term of \eqref{t}. We will look in more detail at the matrix $T$ later in section \ref{T}. Finally, we introduce a measurement matrix $\mathcal{M} \in \mathbb{R}^{N_s \times 39 N}$ that maps the field $X_n$ at each time step to the measured data at the sensors (i.e. implements formulas \eqref{G} and \eqref{M}). Note that $\mathcal{M}$ for the discretised adjoint is defined not the same as for the continuous formulae. We first consider the map from the source $S$ to the measurements.


\begin{definition}\label{Hdef}
The map $\mathbb{H}: \mathbb{R}^{39NN_t} \rightarrow \mathbb{R}^{N_s N_t}$ is defined by
\begin{align} \label{for_map_dis}
\begin{split}
\hat{P}=\mathbb{H} S, \quad \hat{P}_{n}            =\mathcal{M}X_n \quad (n \in \{0,...,N_t-1\}), \quad
\hat{P}              =\Big[\hat{P}_{n+1} \Big]_{n=-1}^{N_t-2},
\end{split}
\end{align}
where $X_n$ is defined by \eqref{t} with initial condition $X_{-1} = 0$, and $\hat{P} \in \mathbb{R}^{N_s N_t}$ is the time series stack of measured data at iterates $n \in \{0,..., N_t-1\}$.
\end{definition}

\noindent In the next lemma we show how to compute the adjoint of $\mathbb{H}$. Note that in fact this lemma applies more generally for the adjoint of any discretised problem taking the form described here.

\begin{lemma}\label{lem2}
The adjoint $\mathbb{H}^*$ of $\mathbb{H}$ defined in definition \ref{Hdef} is given by
\begin{align}  \label{adj}
\big [X^*_n \big ]_{n=-1}^{N_t-2} =\mathbb{H}^* \hat{P} \in \mathbb{R}^{39 N N_t},\quad X^*_{N_t-1} = 0, \quad X^*_{n-1} = T^* X^*_n + \mathcal{M}^* \hat{P}_n, \quad n \in \{0,..., N_t-1\}.
\end{align}
\begin{proof}
Let us assume that $X_n$ satisfies \eqref{t} with initial condition $X_{-1} = 0$, and $X_n^*$ satisfies the second two equations in \eqref{adj}. Then using the conditions $X_{-1} = 0$ and $X_{N_t-1}^* = 0$ we have
\[
\sum_{n=-1}^{N_t-2} (X_{n+1} - X_n ) \cdot \overline{X^*_n} = \sum_{n=0}^{N_t-1} X_n \cdot \overline{(X_{n-1}^* - X_n^*)}.
\]
Then applying \eqref{t} on the left and \eqref{adj} on the right we have
\[
\sum_{n=-1}^{N_t-2} (TX_n - X_n + S_{n+1/2} ) \cdot \overline{X^*_n} =  \sum_{n=0}^{N_t-1} X_n \cdot \overline{(T^* X_{n}^* - X_n^*+ \mathcal{M}^* \hat{P}_n)}.
\]
Rearranging this slightly gives
\[
\sum_{n=-1}^{N_t-2} (TX_n - X_n) \cdot \overline{X^*_n} + \sum_{n=-1}^{N_t-2} S_{n+1/2}  \cdot \overline{X^*_n} = \sum_{n=0}^{N_t-1}  (T X_n - X_n ) \cdot \overline{X_n^*} +  \sum_{n=0}^{N_t-1} (\mathcal{M}X_n) \cdot \overline{ \hat{P}_n}.
\]
Applying again the conditions $X_{-1} = 0$ and $X_{N_t-1}^* = 0$ we see that
\[
\sum_{n=-1}^{N_t-2} S_{n+1/2}  \cdot \overline{X^*_n} = \sum_{n=0}^{N_t-1} (\mathcal{M}X_n) \cdot \overline{ \hat{P}_n}
\]
which is equivalent to $\langle S, \mathbb{H}^* \hat{P} \rangle_{\mathbb{R}^{N N_t}} = \langle \mathbb{H} S, \hat{P} \rangle_{\mathbb{R}^{N_s N_t}}$, and so completes the proof.

\end{proof}
\end{lemma}

\noindent The forward map actually defined in section \ref{f} is
\[
\mathbb{H}\ \mathcal{S},
\]
and so the adjoint is
\[
\mathcal{S}^*\ \mathbb{H}^*.
\]
If we incorporate a time reversal, which amounts to changing $X_{n}^* \mapsto X_{N_t - 2 - n}^*$ in Lemma \ref{lem2}, as well as including $\mathcal{S^*}$, we obtain the following corollary which gives the full method of calculating the adjoint in our case incorporating time reversal. The sum in \eqref{S*H*} is actually just two terms which can also be used to explain \eqref{adj_data} if we commute the summing operation with the computation of $X^*$.

\begin{corollary} \label{adjcor}
$\mathcal{S}^* \mathbb{H}^*$ can be computed as
\begin{equation} \label{S*H*}
\mathcal{S}^* \mathbb{H}^* \hat{P} =\sum_{n=-1}^{N_t-2} \mathcal{S}^*_{n} X^*_{N_t-2-n}
\end{equation}
where $X^*_{n}$ is determined by
\begin{align} \label{col}
X^*_{-1} = 0, \quad X^*_{n+1} = T^* X_n^* + \mathcal{M}^* \hat{P}_{N_t - 2 - n} \ \mbox{for $n \in \{-1, \ ... \ , N_t - 2\}$}.
\end{align}
\end{corollary}

\subsection{The Matrices $T$ and $T^*$} \label{T}

In this section we write the matrices $T$ and $T^*$ explicitly using the forward model presented in section \ref{f} to show how multiplication by each of them may be computed. Considering corollary \ref{adjcor}, we define the adjoint measured data for the discretised adjoint as $\bar{P}_{n+1/2}^\text{adj}= \Delta t \hat{P}_{n+1/2}^\text{adj}$. To start we can write \eqref{num-vel} and \eqref{num-s} in the condensed forms
\begin{align}  \label{matrix-v}
\bar{v}_{n+\frac{1}{2}}
=  A_v   \big[ A_v \  \bar{v}_{n-\frac{1}{2}} + \Phi   \bar{\sigma}_n  \big],
\end{align}
and
\begin{align}   \label{matrix-s}
\bar{\sigma}_{n+1} =  A_\sigma  \big[ A_\sigma   \bar{\sigma}_n  + \Psi_{\text{dis}} \bar{v}_{n+\frac{1}{2}} - \Psi_\text{abs}  \bar{\sigma}_n  \big]+ \theta S_{n+1/2}
\end{align}
where $\theta$ is a sparse matrix that maps $S_{n+1/2}$ to the space of vector $\bar{\sigma}_n$. Also $A_v$, $A_\sigma$, $\Phi$, $\Psi_{\text{dis}}$, and $\Psi_{\text{abs}}$ are matrices that will be described in more detail below although for now we note that $A_v$ and $A_\sigma$ are both diagonal. Based on this we see that $T$ can be written as the following product of matrices in block form
\begin{equation} \label{Tprod}
T =
 \left (
\begin{matrix}
I_v & 0\\
A_\sigma \Psi_{\text{dis}} & A_\sigma (A_\sigma - \Psi_{\text{abs}})
\end{matrix}
\right )
\left (
\begin{matrix}
A_v^2 & A_v \Phi \\
0 & I_\sigma
\end{matrix}
\right )
\end{equation}
where $I_v$ and $I_\sigma$ are the identity matrices. From \eqref{Tprod} we have
\begin{equation} \label{T*prod}
T^* =
 \left (
\begin{matrix}
A_v^2 & 0\\
\Phi^* A_v & I_\sigma
\end{matrix}
\right )
\left (
\begin{matrix}
I_v & \Psi_{\text{dis}}^* A_\sigma \\
0 & (A_\sigma - \Psi_{\text{abs}}^*)A_\sigma
\end{matrix}
\right ).
\end{equation}
Using Corollary \ref{adjcor}, the above equation gives the updates for the adjoint problem as
\begin{align}
\begin{split}
\bar{v}_{n+1/2}                         &= A_v^2 \left(\bar{v}_{n-1/2}+\Psi_\text{dis}^* A_\sigma \bar{\sigma}_{n} \right)\\
\bar{\sigma}_{n+1}                    &= A_\sigma^2 \bar{\sigma}_{n}+ \Phi^* (A_v)^{-1} \bar{v}_{n+1/2}-\Psi_\text{abs}^* A_\sigma \bar{\sigma}_{n}+{\theta^\prime}^*\mathcal{M}^*\bar{P}_{n+1/2}^\text{adj},
\end{split}
\end{align}
where $\theta^\prime$ is a sparse matrix mapping the space of vector $\bar{\sigma}_n$ to the space of vector $X_n$.
Defining $\tilde{\sigma}= A_\sigma \bar{\sigma}$ and $\tilde{v}= A_v \bar{v}$ gives
\begin{align}  \label{main_dis}
\begin{split}
\tilde{v}_{n+1/2}                         &= A_v  \left( A_v \tilde{v}_{n-1/2}+\Psi_\text{dis}^* \tilde{\sigma}_{n} \right)\\
\tilde{\sigma}_{n+1}                    &= A_\sigma \left( A_\sigma \tilde{\sigma}_{n}+ \Phi^*  \tilde{v}_{n+1/2} -\Psi_\text{abs}^* \tilde{\sigma}_{n}\right)+{\theta^\prime}^*\mathcal{M}^*\bar{P}_{n+1/2}^\text{adj}.
\end{split}
\end{align}

Now let us consider the matrices $A_v$, $A_\sigma$, $\Phi$, $\Psi_{\text{dis}}$, and $\Psi_{\text{abs}}$. First we note that $A_v$ and $A_\sigma$ can be found from \eqref{A}. Though the others can be read off from \eqref{num-vel} and \eqref{num-s}, we will write them down explicitly here in order to show how we can explicitly calculate multiplication by their adjoints as required in \eqref{T*prod}.

Toward this goal, let us define the k-space discretised gradient operator
\begin{align}
{\nabla_{(m,i)}^{\pm} }^{p,s}= \frac{\partial_{p,s}}{\partial x_i^{\pm h(m,i)}},
\end{align}
which is defined by \eqref{grf}. Note that the superscript $p,s$ accounts for $c_0^{p,s}$ used in the k-space method (see \eqref{grf}). Here we are considering ${\nabla_{(m,i)}^{\pm} }^{p,s}$ to be a $N \times N$ matrix, and based on \eqref{grf} and the unitarity of the discrete Fourier transform we see that
\begin{equation}\label{gradadj}
\left ({\nabla_{(m,i)}^{\pm} }^{p,s}\right )^* = - {\nabla_{(m,i)}^{\mp} }^{p,s}.
\end{equation}
Using a pseudospectral method without k-space correction, the gradient operator is the same for compressional and shear parts of the fields, and thus this superscript would not be required in that case. We will also need the symmetrisation operator $\mathbb{S}$ acting in the $ij$ indices defined by
\begin{align}
\mathbb{S}[a]_{ij} = \frac{a_{ij} + a_{ji}}{2}.
\end{align}
We apply $\mathbb{S}$ to objects having more indices, but specify that it always acts on the pair $ij$.

From \eqref{num-vel}, the matrix $\Phi \in \mathbb{R}^{9N\times 30N}$ acts on $\tilde{\sigma}$ by
\begin{align} \label{ph1}
{_m(\Phi \tilde{\sigma})}_{i}=\sum_{p,s} {_m(\Phi^\prime \tilde{\sigma})}^{p,s}_{i},
\end{align}
where
\begin{align}  \label{ph2}
{_m(\Phi^\prime \tilde{\sigma})}^{p,s}_{i} = \frac{\Delta t}{\bar{\rho}_i} \left ({\nabla_{(i,m)}^{+} }^{p,s}\right ) \sum_{m'=1}^d {_{m'}\tilde{\sigma}_{im}^{p,s}}.
\end{align}
Here and in what follows we do not include the spatial index $\zeta$ explicitly, but understand that for every value of the other indices ($i$ and $m$ here) we have a vector of dimension $N$, and that the discretised gradient and multiplication by $1/\bar{\rho}_i$ are implemented as operators acting on this spatial index.
Thus, from \eqref{main_dis} and using \eqref{gradadj},
\begin{equation} \label{Phiadj}
{_m(\Phi^{^*} \tilde{v})}^{p,s}_{ij} = -\Delta t\  \mathbb{S} \left [ \left ({\nabla^-_{(i,j)}}^{p,s} \right ) \frac{1}{\bar{\rho}_i} {_j\tilde{v}_i} \right ].
\end{equation}
The symmetrisation $\mathbb{S}$ must be added since this should map into the space of symmetric tensors.

Next, from \eqref{num-s} the operator $\Psi_\text{dis} \in \mathbb{R}^{30N \times 9N}$ acts on $\tilde{v}$  by
\begin{align}
\begin{split}
& {_m(\Psi_{\text{dis}} \tilde{v})}^{p,s}_{ij}  = \Delta t \Bigg [ \delta_{ij}\left [ \bar{\lambda} \left ( {\nabla_{(m,m)}^-}^{p,s}\right )+ \bar{\chi}  \left ( {\nabla_{(m,m)}^-}^{p,s}\right ) \bar{\tau}_{m,\text{dis}}^{p,s} \bar{Y}_{\text{dis}} \right ] Q^{p,s} \sum_{m'=1}^d {_{m'}\tilde{v}_m} \\
&+ 2\ \mathbb{S} \left [ \left [ \bar{\mu}_{ij} \ \delta_{mj}\ \left ( {\nabla_{(i,j)}^-}^{p,s}\right ) + \bar{\eta}_{ij}\ \delta_{mj}  \left ( {\nabla_{(i,j)}^-}^{p,s}\right ) \bar{\tau}_{i,\text{dis}}^{p,s} \ \bar{Y}_{\text{dis}} \right ]Q^{p,s} \sum_{m'=1}^d {_{m'}\tilde{v}_i}\right ] \Bigg ]
\end{split}
\end{align}
From this we can find the formula for the action of the adjoint
\begin{equation} \label{disadj}
\begin{split}
{_m(\Psi_{\text{dis}}^* \tilde{\sigma})}_i & = -\Delta t \sum_{p,s} \sum_{j=1}^d Q^{p,s} \Bigg[ \bigg[   \left ( {\nabla_{(i,i)}^+}^{p,s}\right ) \bar{\lambda} + \bar{Y}_{\text{dis}} \bar{\tau}_{i,\text{dis}}^{p,s} \left ( {\nabla_{(i,i)}^+}^{p,s}\right ) \bar{\chi} \bigg ] {_i\tilde{\sigma}_{jj}^{p,s}} \\
& + 2 \bigg [ \left ( {\nabla_{(i,j)}^+}^{p,s}\right )\bar{\mu}_{ij}  + \bar{Y}_{\text{dis}} \bar{\tau}_{i,\text{dis}}^{p,s}  \left ( {\nabla_{(i,j)}^+}^{p,s}\right ) \bar{\eta}_{ij} \bigg ]  {_j\tilde{\sigma}_{ij}^{p,s}} \Bigg].
\end{split}
\end{equation}
Additionally, from \eqref{num-s}, \eqref{num-s2} and \eqref{ph2}, the operator $\Psi_{\text{abs}} \in \mathbb{R}^{30N \times 30N}$ acts on $\tilde{\sigma}$ by
\begin{align}
\begin{split}
{_m(\Psi_{\text{abs}} \tilde{\sigma})}^{p,s}_{ij} =  {_m(\Psi_\text{abs}^\prime \big( \frac{1}{\Delta t}\big) A_v  \Phi^\prime \tilde{\sigma})}^{p,s}_{ij},
\end{split}
\end{align}
which is actually the action of $\Psi_\text{abs}^\prime$ on $\partial_t \tilde{v}$ using \eqref{num-s} in the form
\begin{align}
\begin{split}
& {_m(\Psi_{\text{abs}}^\prime \partial_t \tilde{v})}^{p,s}_{ij}  =  \Delta t \bigg[ \delta_{ij}\bar{\chi}  \left ( {\nabla_{(m,m)}^-}^{p,s}\right ) \bar{\tau}_{m,\text{abs}}^{p,s} \bar{Y}_{\text{abs}}\partial_t  \sum_{m'=1}^d  {_{m'}\tilde{v}_m^{p,s}} \\
&+ 2\ \mathbb{S} \   \bar{\eta}_{ij}\ \delta_{mj}  \left ( {\nabla_{(i,j)}^-}^{p,s}\right ) \bar{\tau}_{i,\text{abs}}^{p,s} \ \bar{Y}_{\text{abs}}  \partial_t \sum_{m'=1}^d {_{m'}\tilde{v}_i^{p,s}}\bigg]
\end{split}
\end{align}

From \eqref{main_dis}, the action of the adjoint is then given by
\begin{equation} \label{absadj}
\begin{split}
{_m(\Psi_\text{abs}^* \tilde{\sigma})}^{p,s}_{ij}  = _m( \frac{1}{\Delta t} \Phi^{\prime^*}  A_v \Psi_{\text{abs}}^{\prime^*}  \tilde{\sigma} )_{ij}^{p,s},
\end{split}
\end{equation}
where
\begin{align}  \label{absadjp}
\begin{split}
_m({  \Psi_\text{abs}^{\prime^*} \tilde{\sigma})}^{p,s}_{i} & =  - \Delta t \sum_{j=1}^d  \Bigg[   \bar{Y}_{\text{abs}} \bar{\tau}_{i,\text{abs}}^{p,s} \left ( {\nabla_{(i,i)}^+}^{p,s}\right ) \bar{\chi} \ {_i\tilde{\sigma}_{jj}^{p,s}} \\
& + 2     \bar{Y}_{\text{abs}} \bar{\tau}_{i,\text{abs}}^{p,s}  \left ( {\nabla_{(i,j)}^+}^{p,s}\right ) \bar{\eta}_{ij}  \ {_j\tilde{\sigma}_{ij}^{p,s}}\Bigg].
\end{split}
\end{align}
Finally, plugging \eqref{absadj} into the second line in \eqref{main_dis} gives
\begin{align}  \label{main_dis2}
\begin{split}
\tilde{v}_{n+1/2}                 &= A_v  \left( A_v \tilde{v}_{n-1/2}+\Psi_\text{dis}^* \tilde{\sigma}_{n} \right)\\
\tilde{\textbf{v}}_{n+1/2}&=\tilde{v}_{n+1/2}- \frac{1}{\Delta t}  A_v \left(\Psi_{\text{abs}}^{\prime^*} \tilde{\sigma}_n\right)^{p,s} \\
\tilde{\sigma}_{n+1}                 &= A_\sigma \left( A_\sigma \tilde{\sigma}_{n} + \Phi^{\prime^* } \tilde{\textbf{v}}_{n+1/2}\right)+ {\theta^\prime}^*   \mathcal{M}^*\bar{P}_{n+1/2}^\text{adj},
\end{split}
\end{align}
where $\tilde{\textbf{v}}_{n+1/2} \in \mathbb{R}^{2Nd^2}$ is composed of the components $_{m}\textbf{v}_{(i;\zeta;n+1/2)}^{p,s} \in \mathbb{R}^n \ (i,m \in \{1,...,d\})$.
The numerical computation of the above formulae gives the same formulae as in section \ref{nmethodi}. This indicates that using a k-space pseudo-spectral method the numerical computation of the continuous adjoint matches the algebraic adjoint.

\section{First-order Optimization Methods for PAT}    \label{Thirty-six}                                                                                                                                                               We incorporate the forward and adjoint pair in an inverse solver based on the Iterative Shrinkage Thresholding Algorithm (ISTA), which is popular in PAT, e.g. \cite{Arridge,Javaherian}. A fast variant of this algorithm has also been used in PAT \cite{Arridge-b,Huang,Javaherian}.
Let the discretised variant of the sought after initial pressure $P_0$ be denoted by $P \in \mathbb{R}^N$. The inverse problem of inferring $P_0$ from $\hat{P}$ can be fit into a general class of non-smooth constrained convex minimization problems of the form
\begin{align}                                                                                                                                                                                                                          \operatornamewithlimits{argmin}\limits_{P}\left\{F(P):=f(P)+g(P)\right\},   \label{two}                                                                                                                                                \end{align}                                                                                                                                                                                                                            where                                                                                                                                                                                                                                  
$\providecommand{\norm}[1]{\left\lVert#1\right\rVert}                                                                                                                                                                                  f(P)= \frac{1}{2}\norm{\mathbb{H}P-\hat{P}}^2$ is a continuously differentiable function with Lipschitz continuous gradient having smallest Lipschitz constant $\providecommand{\norm}[1]{\left\lVert#1\right\rVert}
L_f={\mathcal{S}_i}_{\text{max}}(\mathbb{H}^T\mathbb{H})$ with ${\mathcal{S}_i}_{\text{max}}(.)$ the largest singular value. The gradient of $f$ is computed by
\begin{align}
\nabla{f}(P)=\mathbb{H}^*\left(\mathbb{H}P-\hat{P}\right).    \label{Six}                                                                                                                                                                                \end{align}                                                                                                                                                                                                                           
Using a total variation (TV) regularized variant of ISTA, we take                                                                                                                                                                                                     $g(P)= \lambda_r \mathcal{J}(P)+\delta_C\left(P\right)$, where $\mathcal{J}(P)$ represents a TV penalty functional, $\lambda_r$ denotes the regularization parameter, and $\delta_C$ is an indicator function for the set of constraints $C = \{P \geqslant 0\}$ \cite{Goldstein,Beck}.

Applying the so-called forward-backward splitting method to a fixed point iterative scheme arising from the optimality conditions of problem \eqref{two} gives two-steps at each iteration $k$ of the optimization algorithm. The first step uses a steepest descent search direction $- \nabla f(P^{k-1})$ and step size $\Gamma^k$
in the form
\begin{align} \label{smooth}
Y^k = P^{k-1}-\Gamma^k \ \nabla{f}(P^{k-1}),                                                                                                                                                                                            \end{align}                                                                                                                                                                                                                            and is called the forward gradient descent step \cite{Goldstein,Beck}. Applying ISTA, the iterates $P^k$ are converged to a minimizer $P^*$ of
problem \eqref{two} if $\Gamma^k \in \left(0, 2/L_f\right)$ \cite{Beck}. Here, $L_f$, the largest singular value of $\mathbb{H}^T\mathbb{H}$, is computed iteratively by the power method following \cite{Arridge,Arridge-b,Javaherian}. Since $L_f$ is agnostic to the unknown $P_0$, it can be stored and used for other experiments done in a fixed setting \cite{Arridge-b}.
The second step is a proximal map in the form
\begin{align} \label{proxmal}                                                                                                                                                                                                          \providecommand{\norm}[1]{\left\lVert#1\right\rVert}                                                                                                                                                                                   \begin{split}                                                                                                                                                                                                                          \text{prox}_{\Gamma^k}(g)(Y^k):=                                                                                                                                                                                                  \operatornamewithlimits{argmin}\limits_{P \geqslant 0}\left\{ g(P)+\frac{1}{2\Gamma^k}\norm{P-Y^k}^2 \right\},
\end{split}                                                                                                                                                                                                                            \end{align}                                                                                                                                                                                                                            
and is called backward gradient step \cite{Goldstein}. Following \cite{Huang,Arridge}, here the proximal map is computed based on Chambolle's dual approach (See \cite{Beck}).
In our study, we will terminate ISTA if the following criterion is satisfied:
\begin{align}
k>1 \ \cap \ 1-\frac{F^k}{F^{k-1}}< \epsilon.
\end{align}
Here, $\epsilon$ is a stopping tolerance, and is chosen close to zero.

\section{Numerical results}
\noindent
The numerical implementation of the system of coupled first-order equations that describe the propagation of PA waves in linear isotropic elastic media based on the pseudo-spectral time-domain method is available on the \textit{k-Wave} website \cite{treeby-b,Treeby-tool,Firouzi}. We modified this code so that it includes the absorption and physical dispersion following the frequency power law, using the splitting technique, as discussed in section \ref{f} \cite{Treeby-f}. To validate the computed forward and adjoint models, \eqref{adj_main} was first used to check if the inner product test is satisfied for any initial pressure $P_0$ and data. We then performed reconstructions from simulated data in both 2D and 3D settings as described below.

\subsection{2D phantom} \label{2d}

\subsubsection{Computational grid}
we used a computational grid with a size of $14 \times 14 \ \text{cm}^2$ to simulate the size of the top surface of the skull.

\noindent
\textit{Data generation:} To simulate the propagation of wavefields, the computational grid was made up of $472\times 472$ grid points equidistantly spaced with a
separation distance of $2.9661 \times 10^{-2} \ \text{cm}$ along both Cartesian coordinates. This computational grid was enclosed by a PML having a thickness of 20 grid points and a maximum attenuation coefficient of 2 nepers per grid point so that a good trade-off between mitigating spurious wave wrapping at the boundaries and reflection of waves at the edge of the PML was made \cite{Tabei}. The propagated pressure field was measured in time by $200$ detectors that were evenly placed aligned by the top half of periphery of a circle having a radius of $r=6.8\text{cm}$ so that $\pi$ radians were covered by the detectors. The skull was simulated with semi-circular interfaces with distances of $0.85r$ and $0.95r$ to the center of the semi-circle so that it has an even thickness of $6.8\text{mm}$. This has provided an even radial distance of 3.4mm between the outer edge of the skull and the detectors.

\noindent
\textit{Image reconstruction:}
To avoid an inverse crime for discretisation \cite{Kaipio}, the image reconstruction was done on a grid made up of $328 \times 328$ grid points which are placed evenly with a separation distance of $4.2683 \times 10^{-2} \ \text{cm}$ along both Cartesian coordinates. The thickness of the PML was reduced to 16 grid points.

\noindent
\subsubsection{Physical parameters}

\noindent
The maps corresponding to the medium's mass density $\rho$, compressional wave propagation speed $c_p$ and shear wave propagation speed $c_s$ were shown in figures \ref{F1a}, \ref{F1b} and \ref{F1c}, respectively. The colour scales are shown to the right of each map, where the blue colour represents the physical parameters of soft tissue with $c_p=1500 \ \text{ms}^{-1}$, $c_s=0 \ \text{ms}^{-1}$ and $\rho=1000 \ \text{kgm}^{-3}$, and the red colour represents the skull with $c_p=3000 \ \text{ms}^{-1}$, $c_s=1500 \ \text{ms}^{-1}$ and $\rho=1850 \ \text{kgm}^{-3}$. These parameters were chosen following \cite{Mitsuhashi}.
The absorption coefficients were set to $\alpha_{0,p}=10 \ \text{dB} \text{MHz}^{-y} \text{cm}^{-1}$ and $\alpha_{0,s}=20 \  \text{dB} \text{MHz}^{-y} \text{cm}^{-1}$ in the skull, and $\alpha_{0,p}=0.75 \
\text{dB} \text{MHz}^{-y} \text{cm}^{-1}$ and $\alpha_{0,s}=0.5 \ \text{dB} \text{MHz}^{-y} \text{cm}^{-1}$ in the soft tissue. Note that we assumed absorption coefficients associated with the skull greater than the experimental values obtained in \cite{White} (cf. Table 1 in \cite{Treeby-f}). Following \cite{Treeby-f}, the exponent factor was assumed constant across the entire medium, and was set to $y=1.4$.

\subsubsection{Validation of adjoint}
Using the setting described above, we numerically measured the accuracy of the computed adjoint model using the inner product test in \eqref{adj_main}. To do this, we used a randomly selected vector for $\hat{P}$, together with an initial pressure distribution $P_0$ in the form of a circular disk with a radius of $0.8r$, where the values at each point of the disk are chosen randomly. The relative difference between the left-hand and right-hand sides of \eqref{adj_main} was averaged between 10 attempts. This gives values $7.43 \times 10^{-5}$ and $8.71 \times 10^{-6}$ for the grids used for image reconstruction and data generation, respectively. Our observations showed us that with an increase in density of the grid, the inner product test is satisfied with a higher order of accuracy.

\subsubsection{Simulation setting}
To evaluate the performance of the forward and adjoint pair for image reconstruction, we considered two scenarios as follows.

\begin{figure}\centering
{\subfigure[]{\includegraphics[width=.325\textwidth]{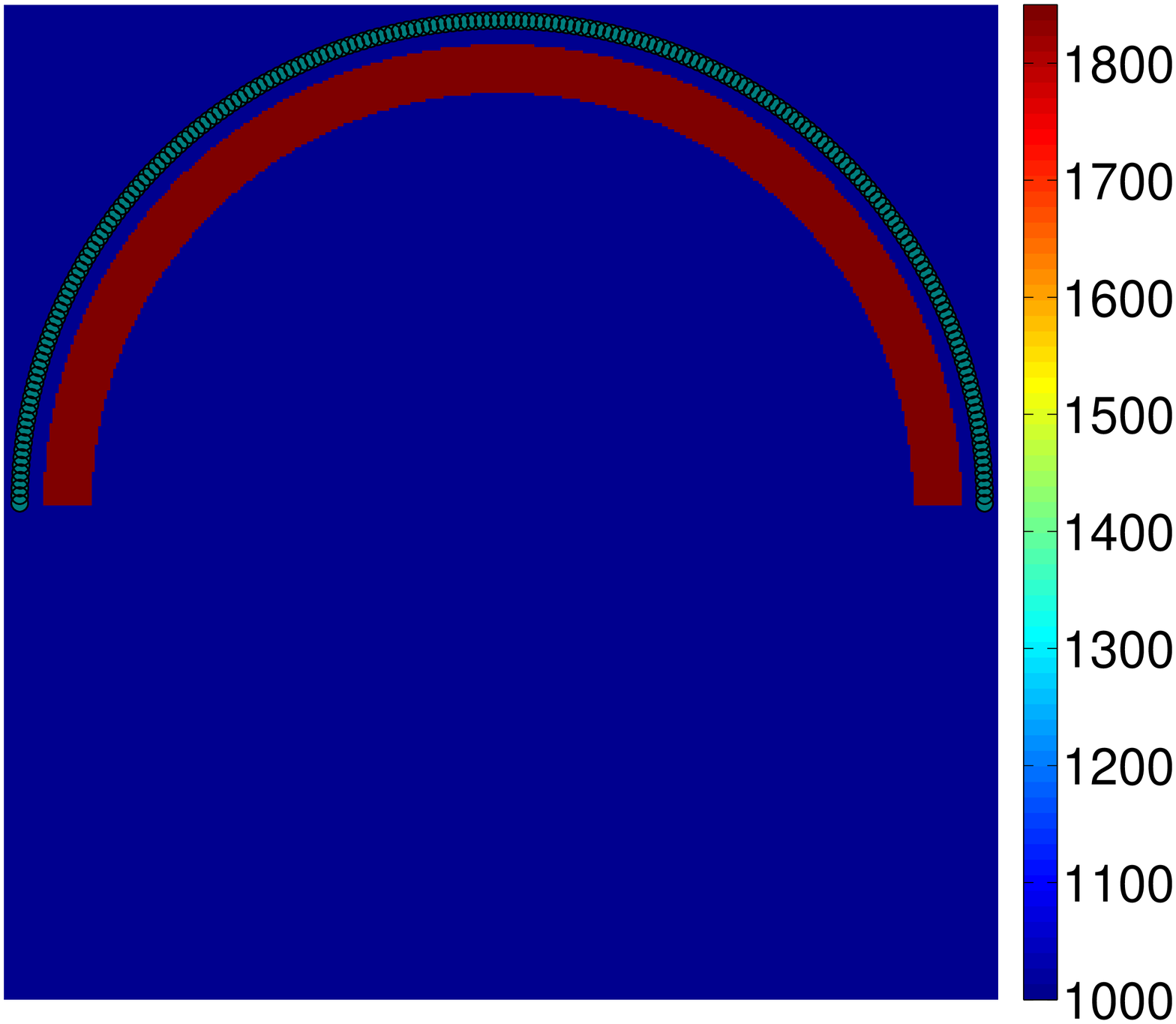}\label{F1a}}
\subfigure[]{\includegraphics[width=.325\textwidth]{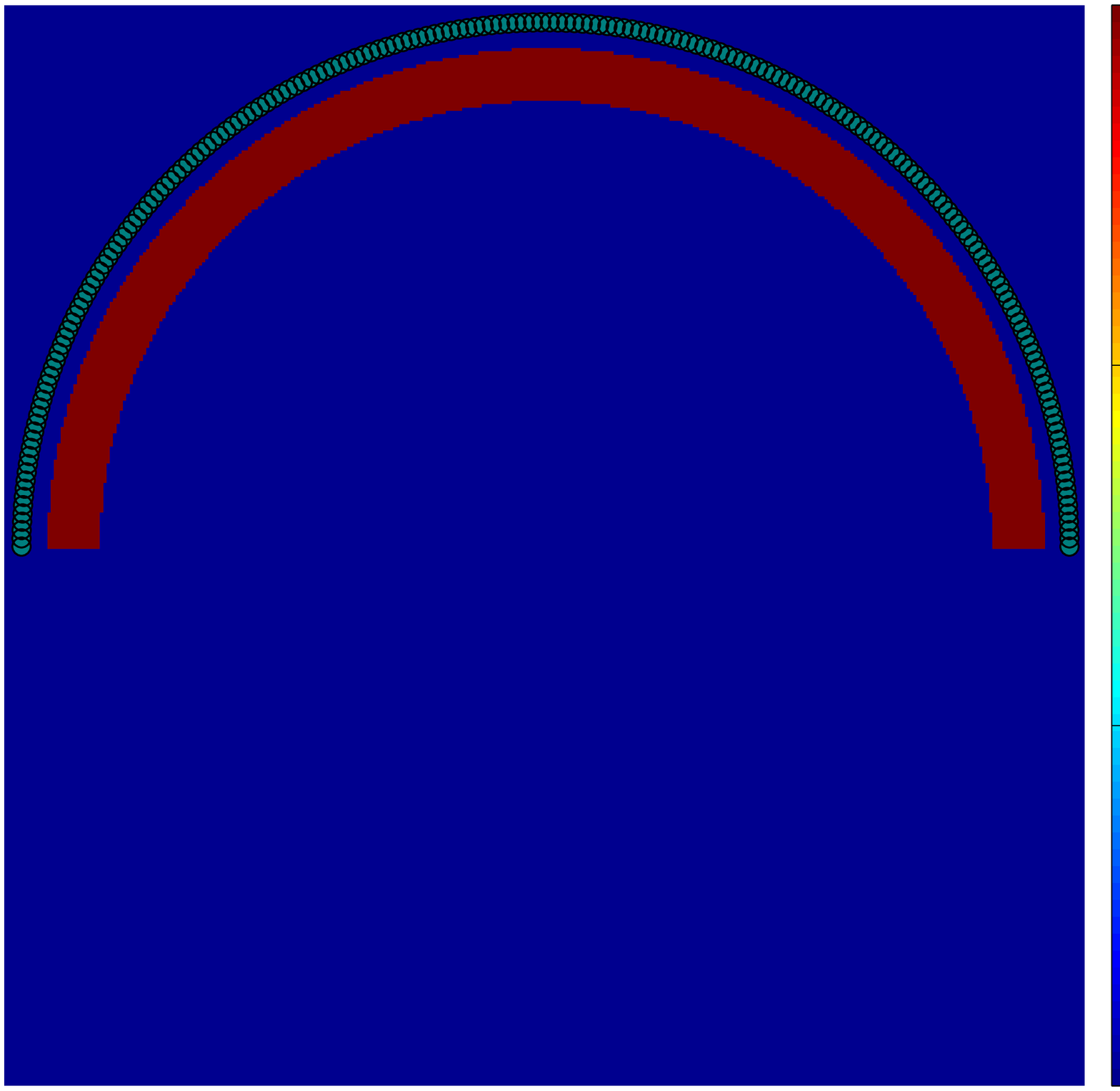}\label{F1b}}
\subfigure[]{\includegraphics[width=.325\textwidth]{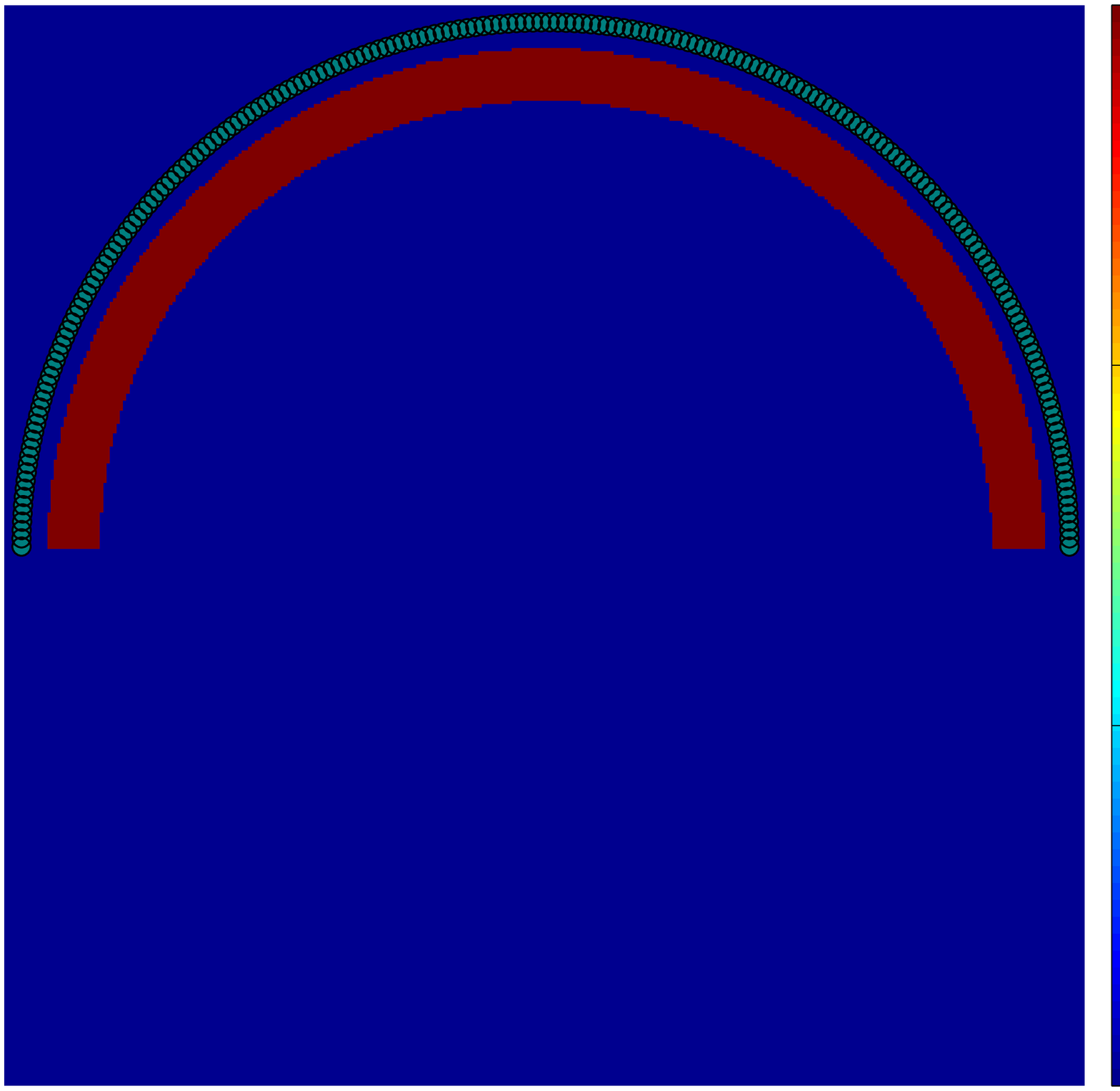}\label{F1c}}\\
\subfigure[]{\includegraphics[width=.325\textwidth]{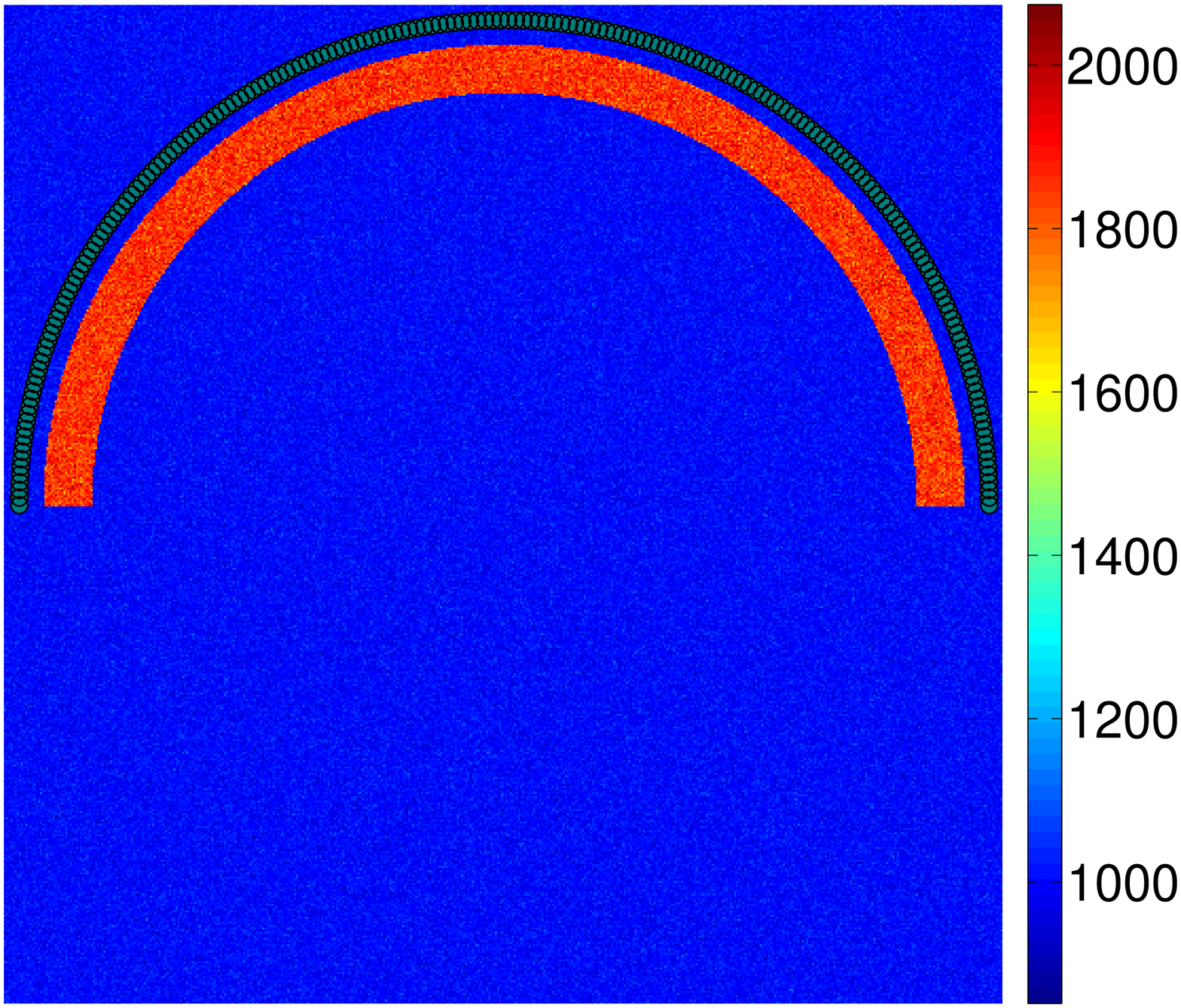}\label{F1d}}
\subfigure[]{\includegraphics[width=.325\textwidth]{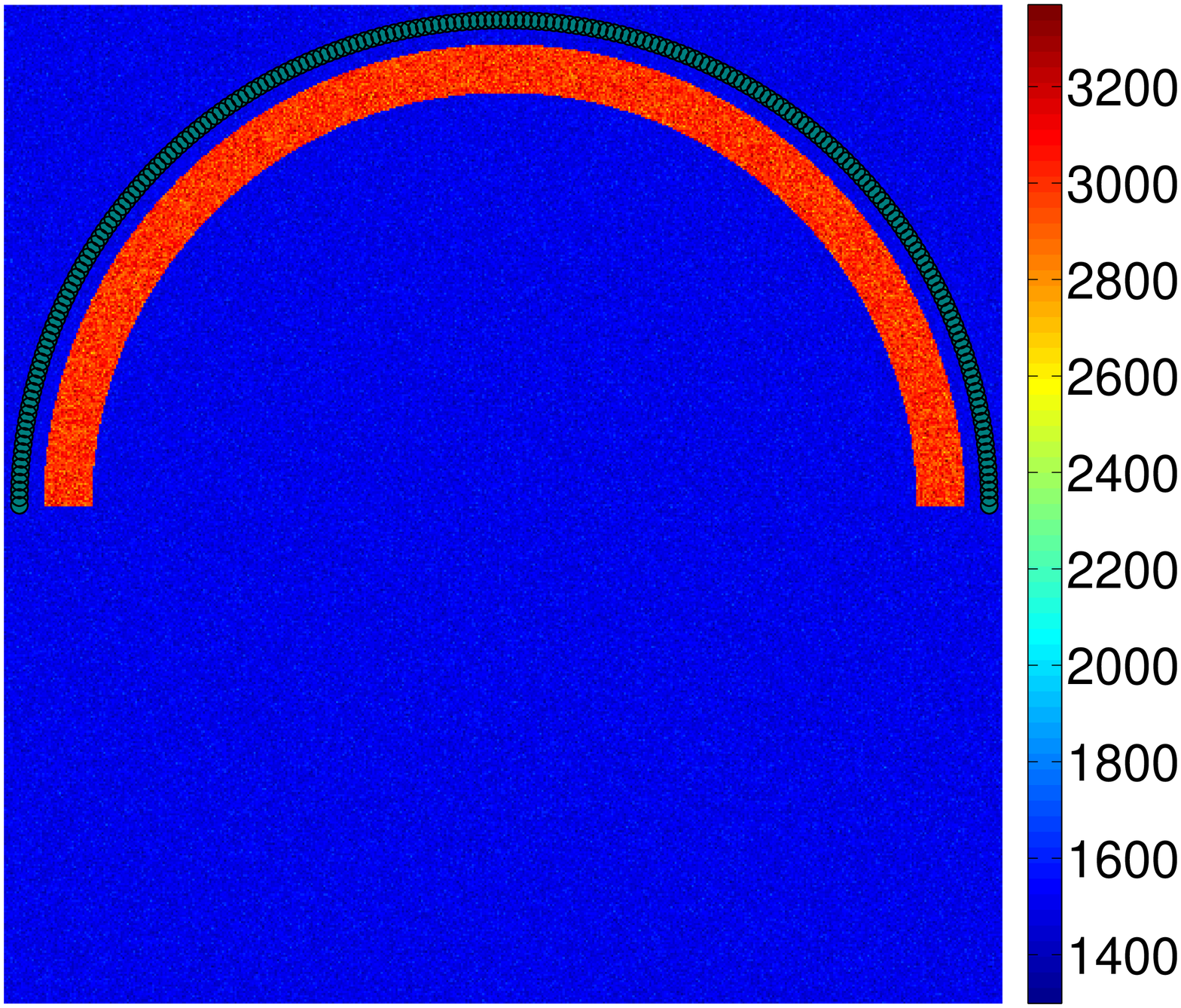}\label{F1e}}
\subfigure[]{\includegraphics[width=.325\textwidth]{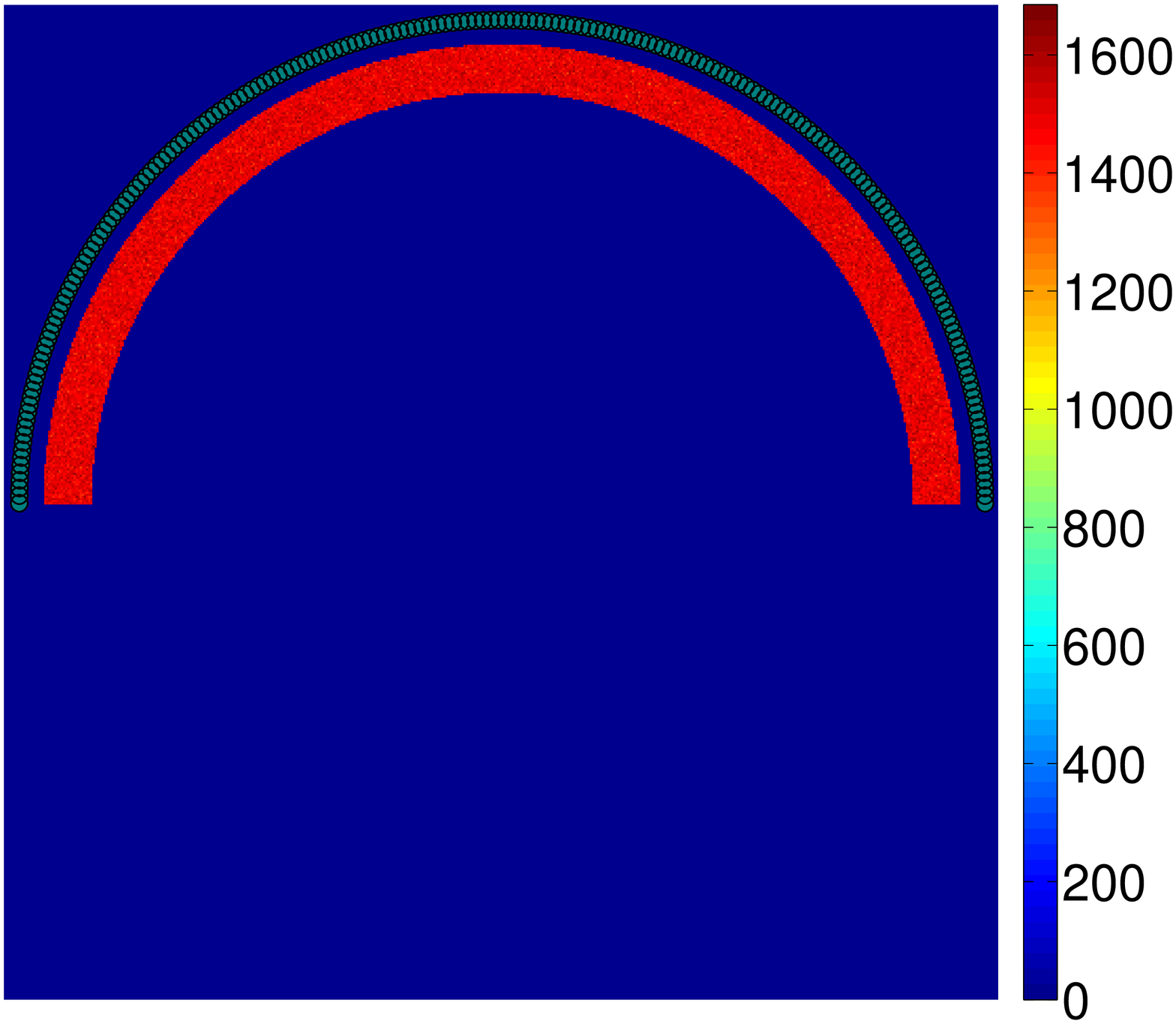}\label{F1f}}
}
\caption{2D phantom. Exact physical maps: (a) $\rho$ (b) $c_p$ (b) $c_s$, and noise-contaminated physical maps: (d) $\rho$ (e) $c_p$ (f) $c_s$.}
\end{figure}

\textit{Scenario1:} In general, the inverse problem in PAT is based on the assumption that the physical parameters of the medium are known. In our first experiment, we used the maps in figures \ref{F1a}, \ref{F1b} and \ref{F1c} as physical parameters for both data generation and image reconstruction. This implies that we have an exact knowledge of the physical parameters. Since this assumption does not hold in practical cases, this is considered as an inverse crime \cite{Kaipio}. Using these maps, the grid used for data generation supports a maximal frequency up to $2.5286 \ \text{MHz}$ for propagation of compressional waves through the entire medium and shear waves within the skull.

\textit{Scenario 2:} In the second experiment, we avoided an inverse crime in estimating medium's parameters by using different maps for data generation and image reconstruction. Correspondingly, for generating data we contaminated the maps in figures \ref{F1a}, \ref{F1b} and \ref{F1c} with a $30\text{dB}$ Additive White
Gaussian Noise (AWGN). The contaminated maps are displayed in figures \ref{F1d}, \ref{F1e} and \ref{F1f}. For image reconstruction, we assumed the contaminated maps are not readily available, and thus we used the clean maps.
Using the noise contaminated maps for data generation, the associated grid supports maximal frequency up to $2.2047 \ \text{MHz}$ for propagation of shear waves within the skull and $2.1889 \ \text{MHz}$ for compressional waves through the entire medium. In these figures, the location of ultrasound detectors has been shown by the green semi-circle.

The grid used for image reconstruction supports a maximal frequency of $1.7571 \ \text{MHz}$ for compressional waves through the entire medium and shear waves within the skull.
We created the initial pressure map with a maximal amplitude of $2$, as shown in figure \ref{F2a}.
For both scenarios, a $CFL$ of 0.3 was sufficient to guarantee the stability of the forward and adjoint models. Accordingly, the computed pressure wavefield was recorded in $4451$ time steps, and interpolated to the detectors using trilinear interpolation \cite{Mitsuhashi}. The generated data (for both scenarios) were then contaminated with a $30$dB AWGN.
\subsubsection{Image reconstruction}
\begin{figure}\centering                                                                                                                                                                                                               {\subfigure[]{\includegraphics[width=.45\textwidth]{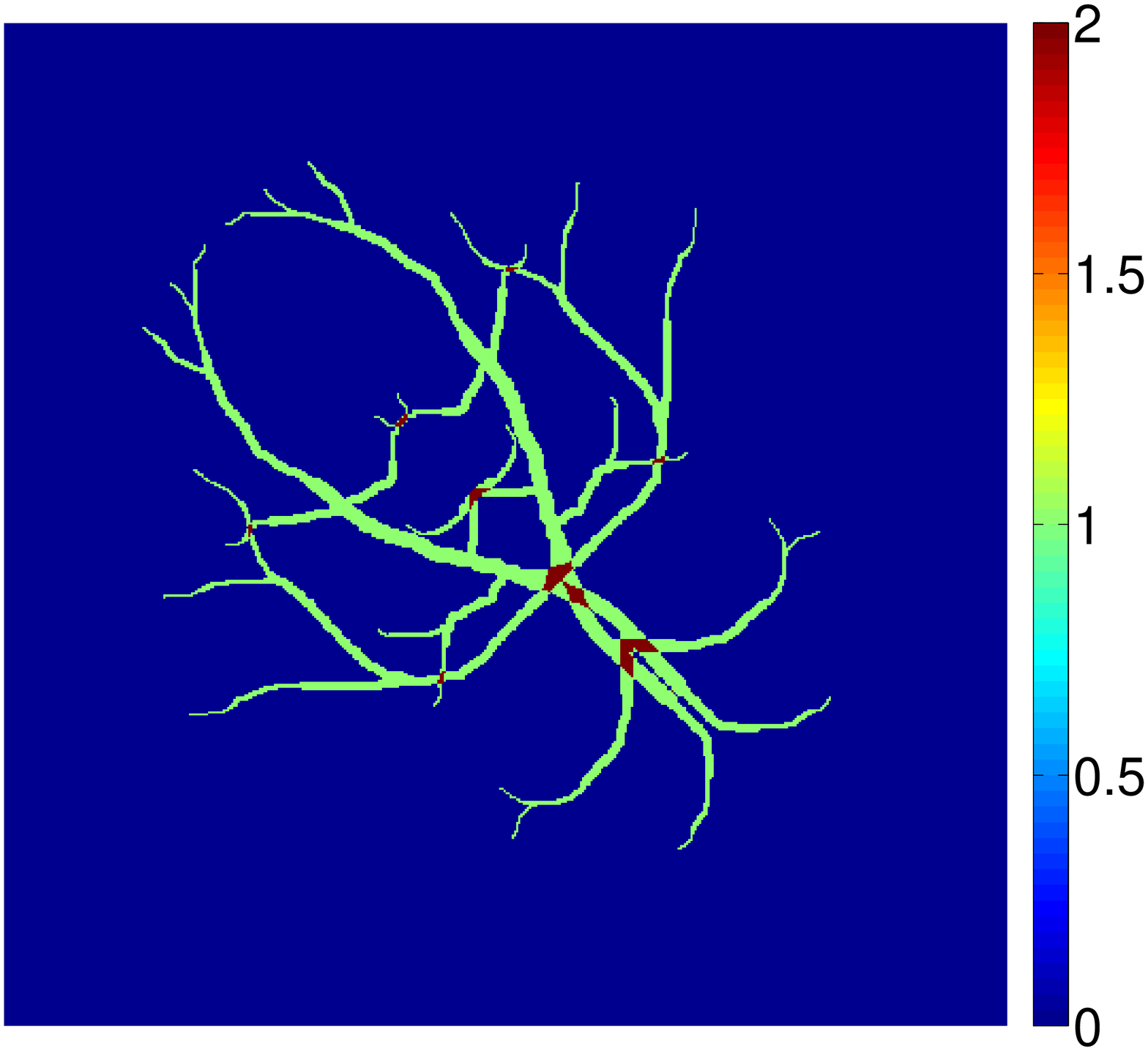}\label{F2a}}\\
\subfigure[]{\includegraphics[width=.45\textwidth]{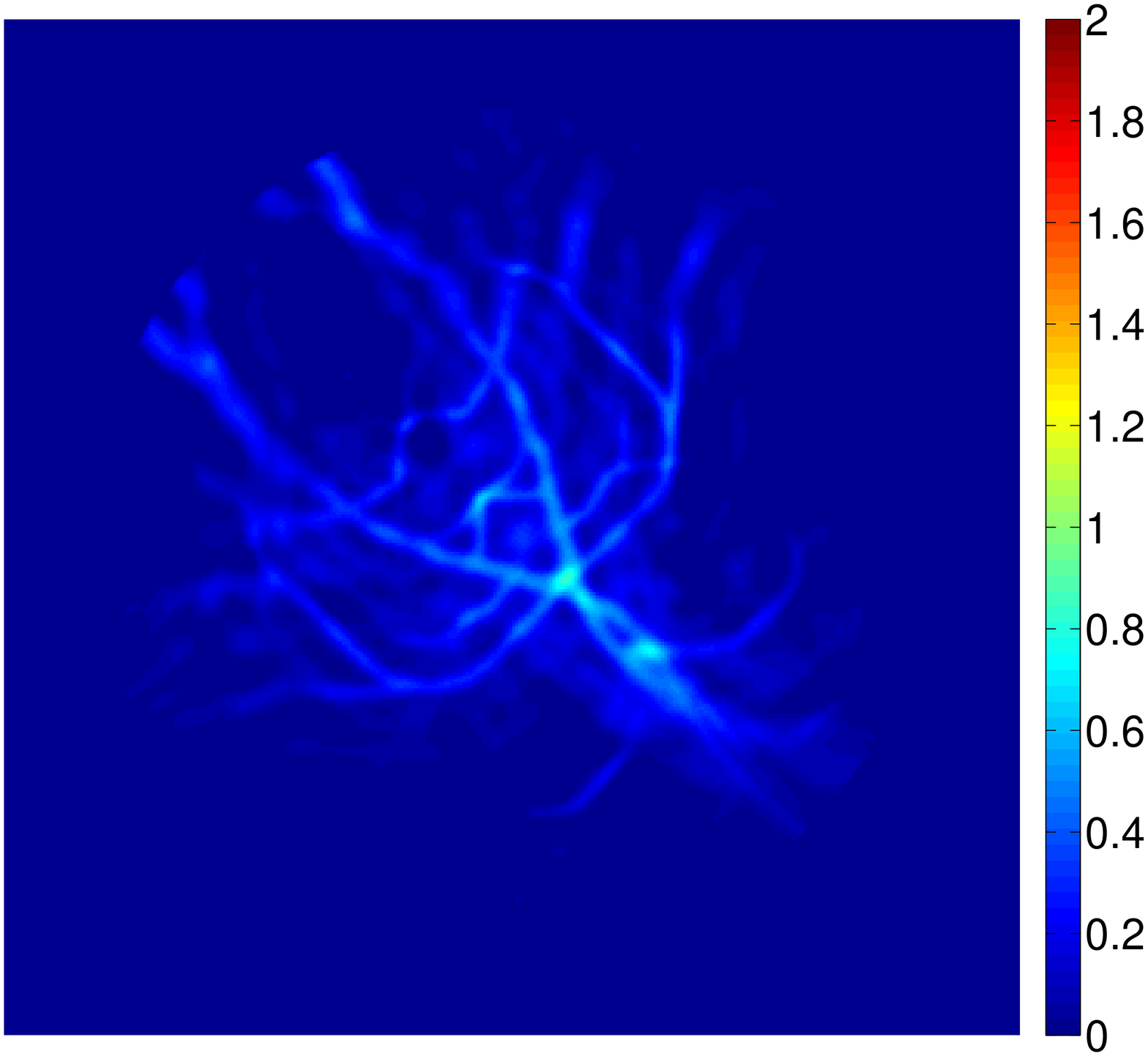}\label{F2b}}
\subfigure[]{\includegraphics[width=.45\textwidth]{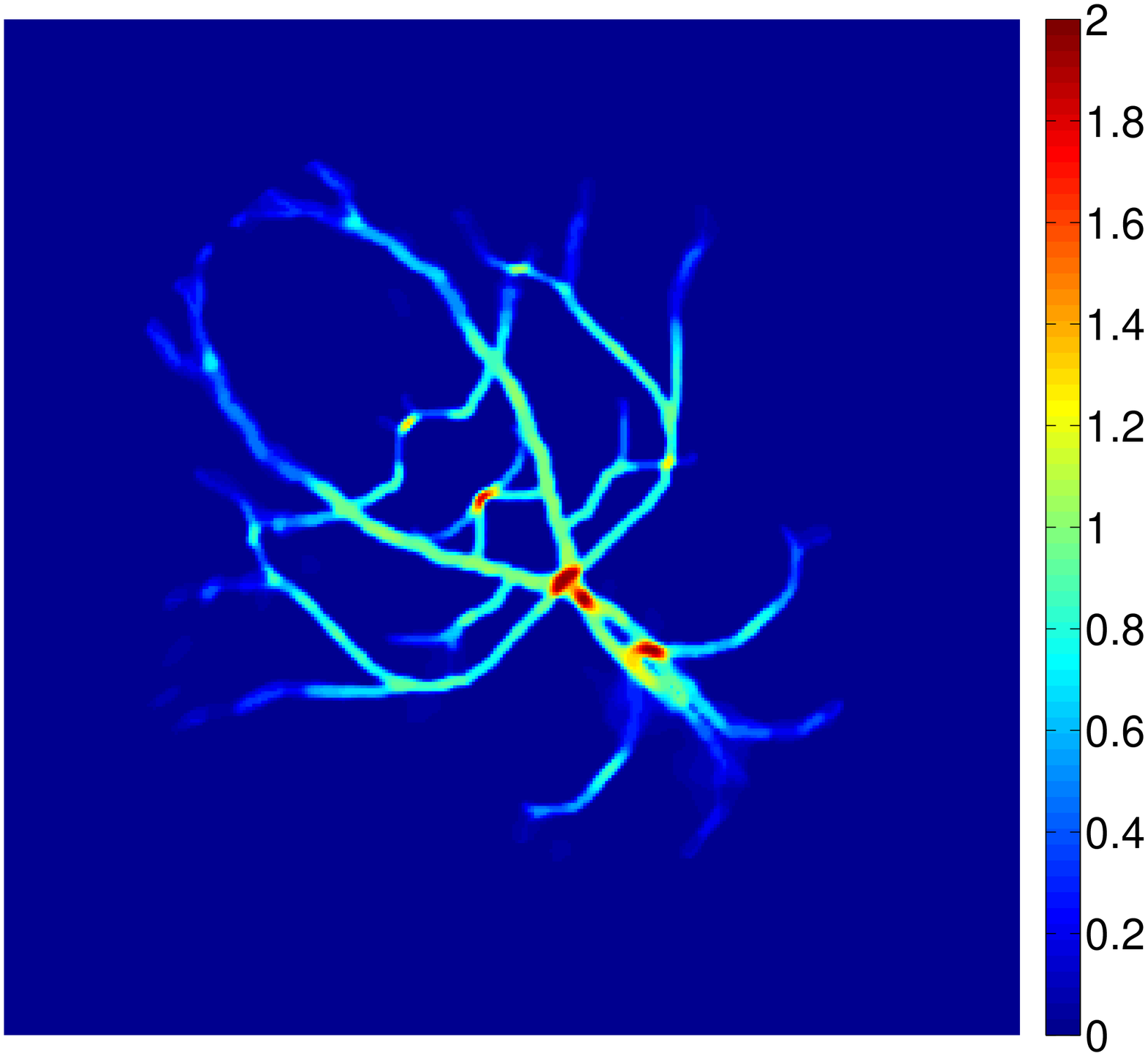}\label{F2c}}\\
\subfigure[]{\includegraphics[width=.45\textwidth]{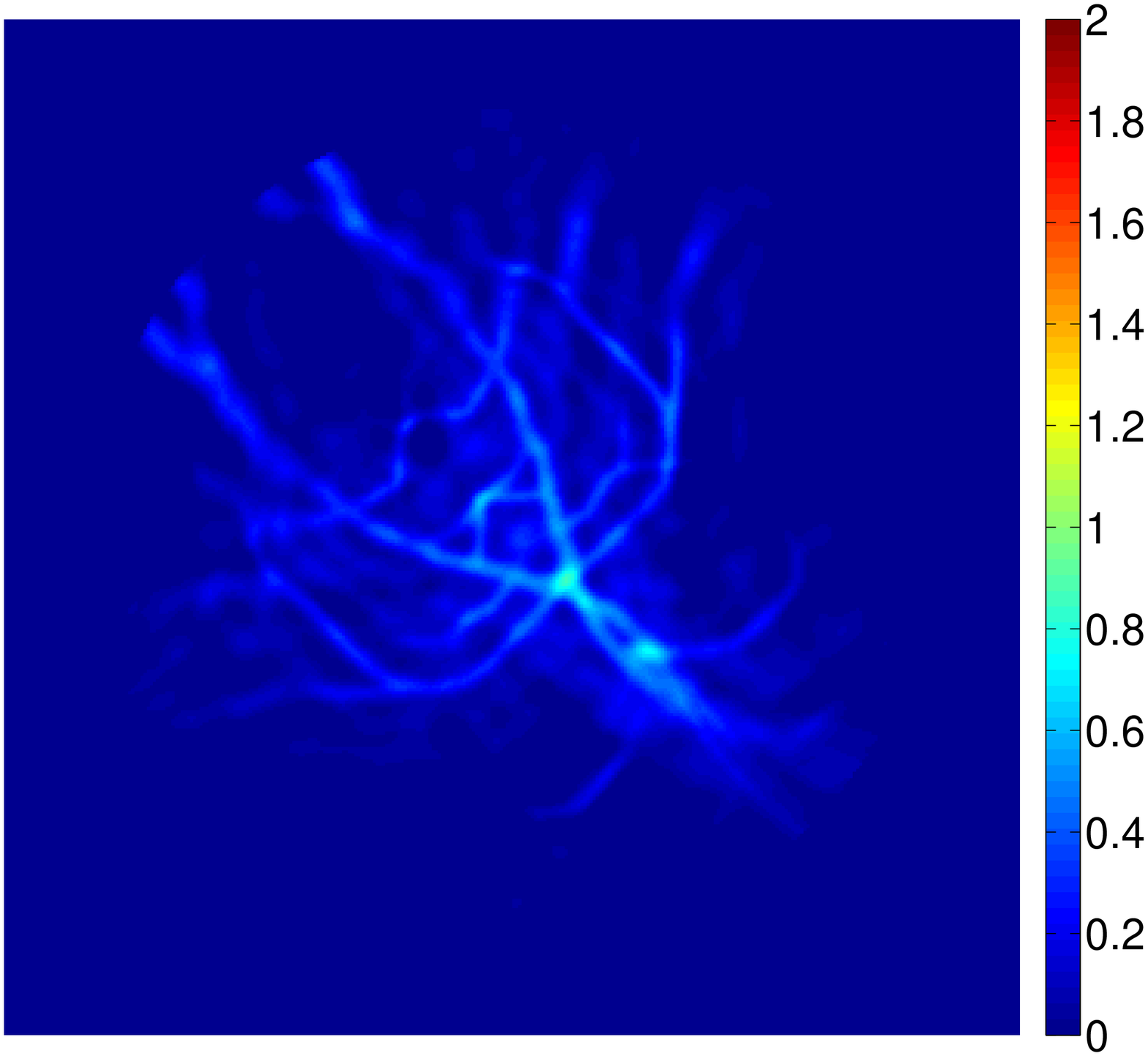}\label{F2d}}
\subfigure[]{\includegraphics[width=.45\textwidth]{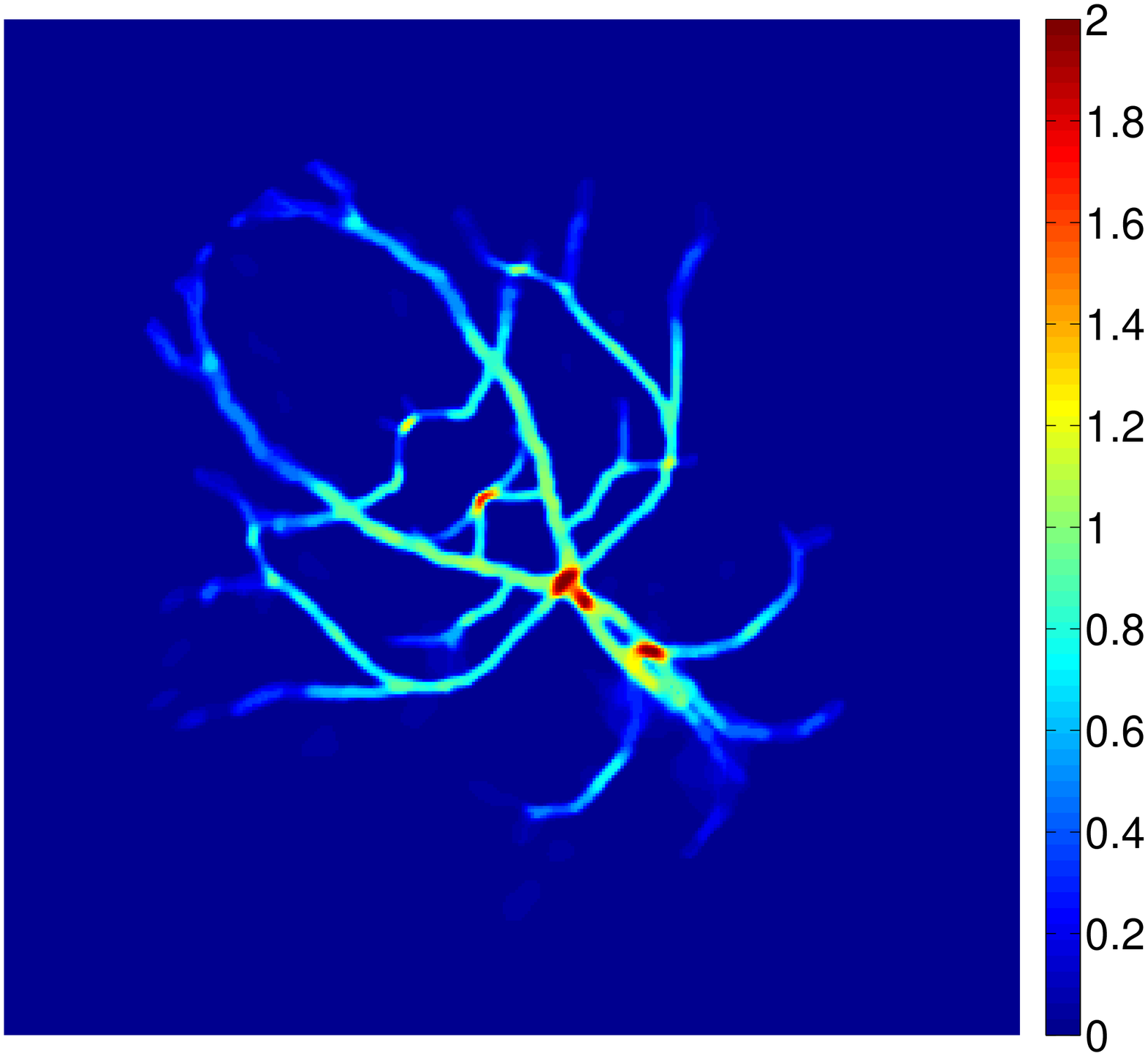}\label{F2e}}
}
\caption{2D phantom. (a) initial pressure map, and reconstructed images using exact physical parameters (inverse crime): (b) TR ($\alpha_{0,p,s}=0$) (c) ISTA, and erroneous physical parameters: (d) TR ($\alpha_{0,p,s}=0$) (e) ISTA.}
\end{figure}

We first reconstructed an image corresponding to each scenario using the time reversal method. This was performed using \textit{k-Wave} toolbox \cite{treeby-b,Treeby-tool}. According to \cite{Treeby}, a filtering of the absorption and dispersion terms in the spatial frequency domain may be required to ensure the stability of TR. Here, to make a fair comparison between TR and ISTA, we applied TR optimistically on a non-absorbing medium with $\alpha_{0,p,s}=0$.
The images reconstructed by TR  for scenarios 1 and 2 are displayed in figures \ref{F2b} and \ref{F2d}, respectively. The computed forward operator $\mathbb{H}$ and adjoint $\mathbb{H}^*$ were then incorporated into the inverse solver discussed in section \ref{Thirty-six}. The regularization parameter was empirically set to $\lambda_r = 1 \times 10^{-2}$. A step size of $\Gamma_k =1.8/L_f$ was chosen and used for all iterations $k$. Here, $L_f$ was computed by the power iteration method \cite{Arridge,Arridge-b,Javaherian}. The iterates of power iteration algorithm converged to $L_f$ after around 15 iterations. For ISTA, the iterates are initialized by zero, and the algorithm was terminated using the stopping tolerance $\epsilon=1 \times 10^{-4}$. Figures \ref{F2c} and \ref{F2e} show an image of the final iterate computed by ISTA for scenarios 1 and 2, respectively.

\noindent
The computed sequence of iterates was measured by two parameters:

\noindent
(1) Relative Error (RE):
\begin{align}                                                                                                                                                                                                                          \providecommand{\norm}[1]{\left\lVert#1\right\rVert}                                                                                                                                                                                   RE(P^k)=\frac{\norm{P^k-{\tilde{P}}_\text{phantom}}_2}{\norm{{\tilde{P}}_\text{phantom}}_2} \times 100,                                                                                                                                                      \end{align}
where $P^k$ denotes the update at iteration $k$, and $\tilde{P}_\text{phantom}$ denotes the initial pressure distribution in the phantom interpolated to the grid used for image reconstruction.

\noindent
(2) Objective function ($F(P^k)$): (cf. section \ref{Thirty-six}, \eqref{two}).

\subsubsection{Observations}
\begin{figure}\centering
{\subfigure[]{\includegraphics[width=.45\textwidth]{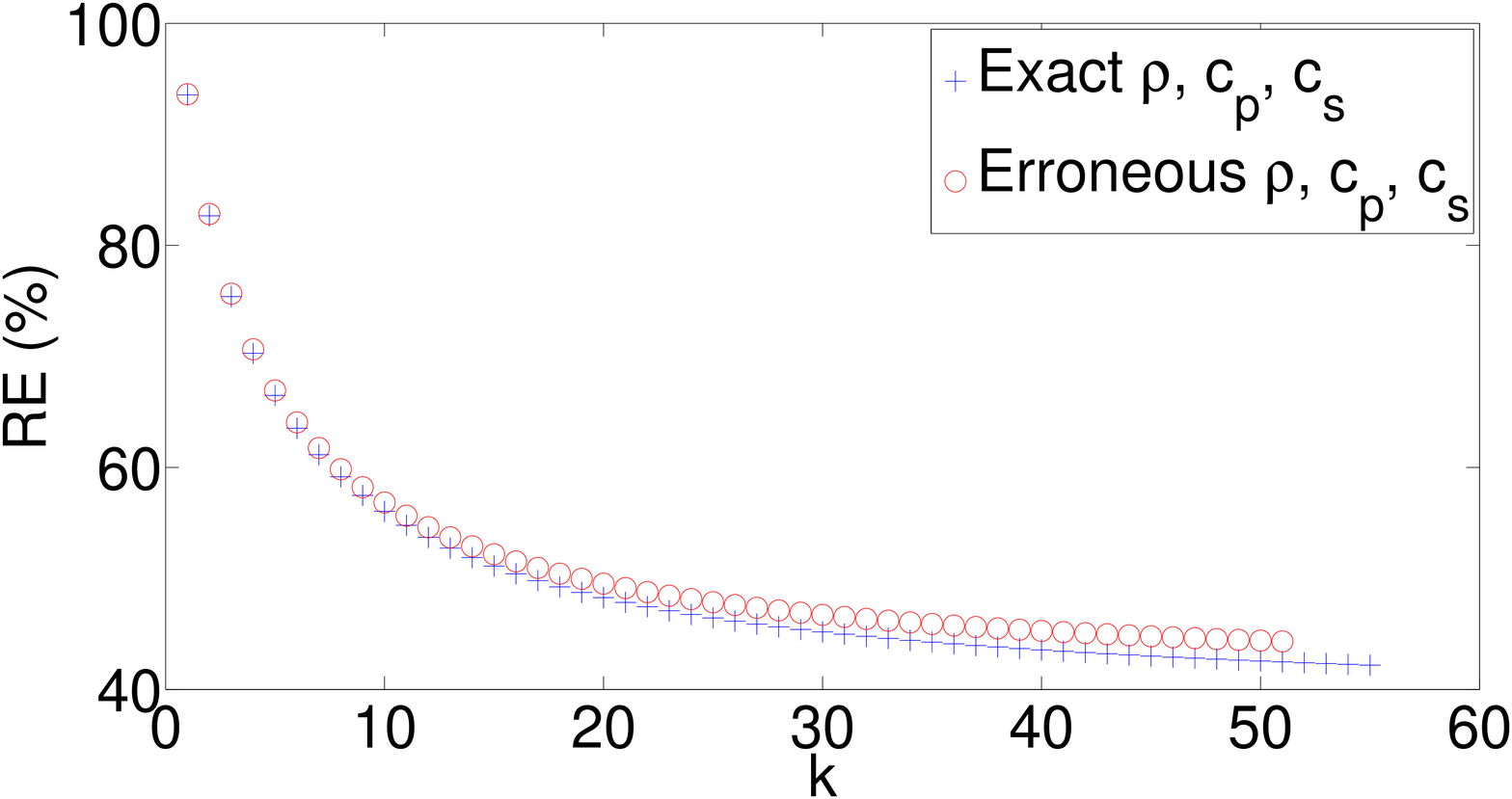}\label{F3a}}\\
\subfigure[]{\includegraphics[width=.45\textwidth]{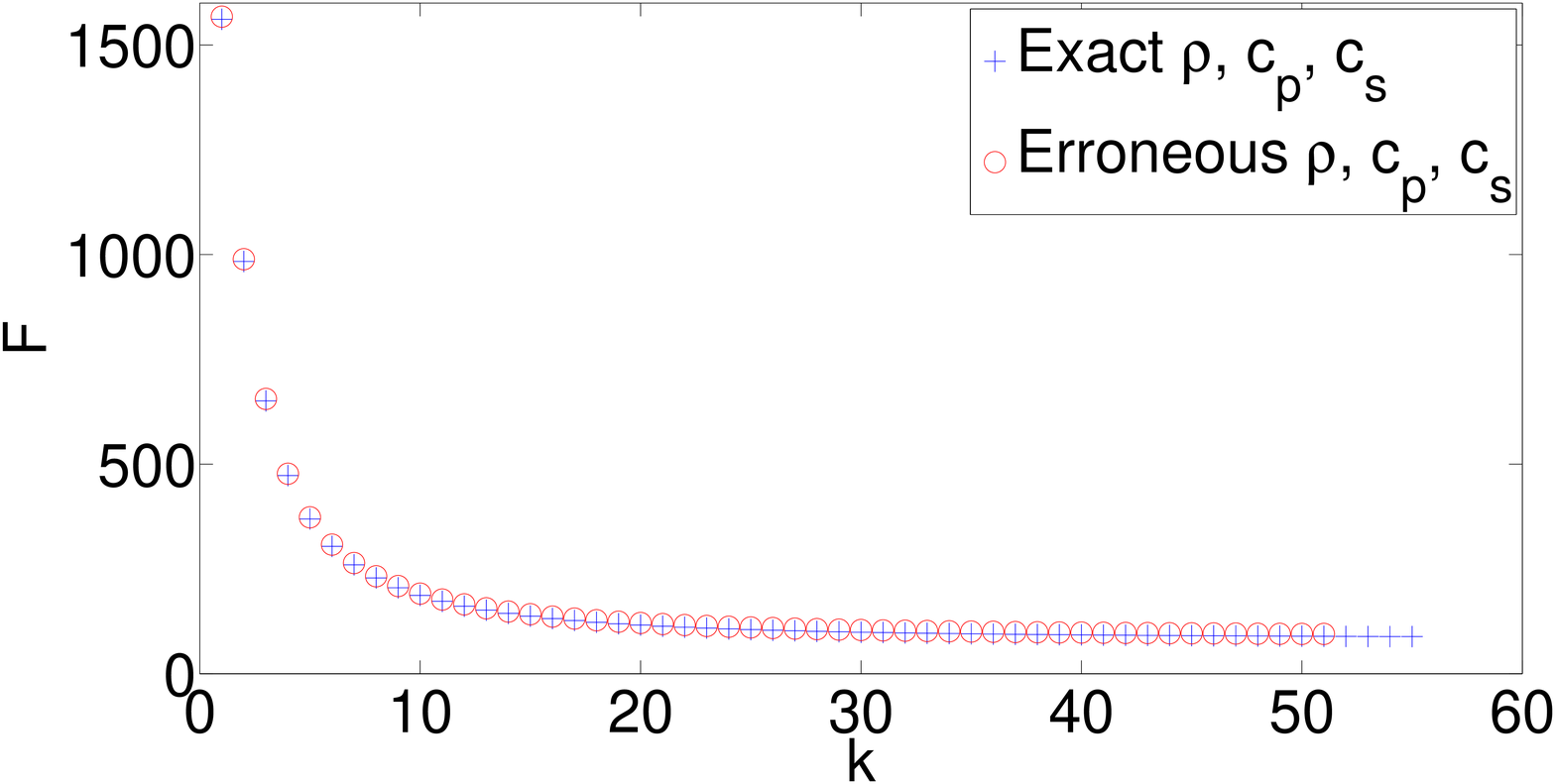}\label{F3b}}
\subfigure[]{\includegraphics[width=.45\textwidth]{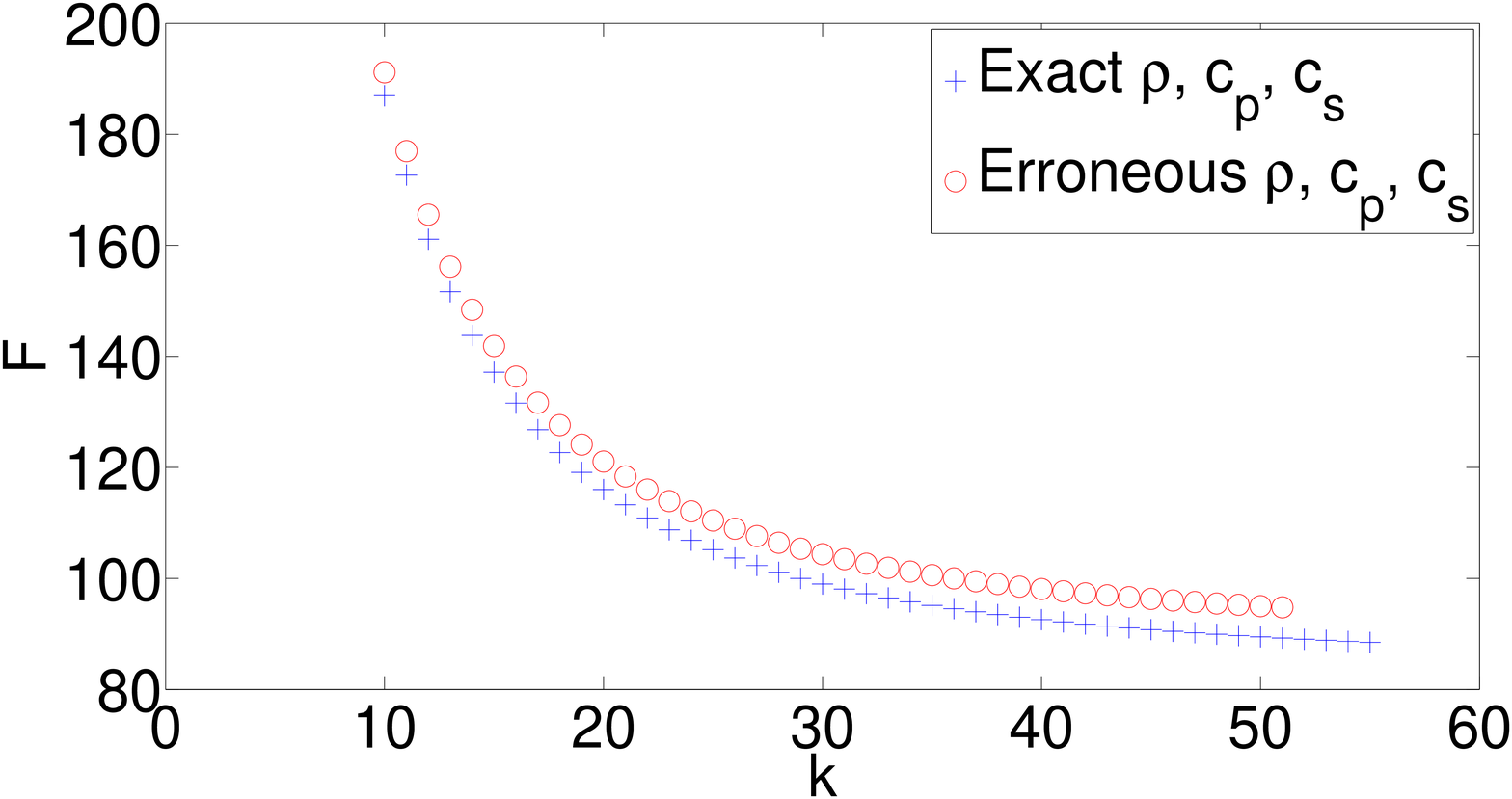}\label{F3c}}
}
\caption{2D phantom. (a) Relative Error (RE) (b) objective function (F) (c) F around the stopping point.}
\end{figure}

Figures \ref{F3a} and \ref{F3b} show RE and objective function values of iterates computed by ISTA versus the iteration number $k$, respectively. Figure \ref{F3c} shows $F$ from a large view around the stopping point. In these figures, the blue and red plots, respectively correspond to scenarios 1 and 2. Our numerical observations for the two mentioned scenarios are as follows.

\noindent
\textit{Scenario 1:}
Both $RE$ and $F$ were monotonically reduced, and the stopping criterion was satisfied at iteration 55. The $RE$ and $F$ reached values of $42.19\%$ and $87.24$ at the final iteration, respectively. The final iteration pertains to the image shown in figure \ref{F2c}. From figures \ref{F3b} and \ref{F3c}, ISTA has reduced $F$ almost $95\%$.

\noindent
\textit{Scenario 2:}
In scenario 2 when we avoided the inverse crime in estimating physical parameters, a monotonic reduction in both $RE$ and $F$ was observed, and the stopping criterion was satisfied at iteration 51. As shown in figures \ref{F3a} and \ref{F3b}, RE and $F$ reached values of $44.29\%$ and $93.52$ at the final iteration, which corresponds to the image shown in figure \ref{F2e}.

These figures indicate that in presence of an error in estimating physical parameters, the inverse solver was tolerant enough to reconstruct almost the same image as using the exact physical maps.

\subsection{3D phantom}

\subsubsection{Computational grid}
The grid was created as a rectangular cuboid with a size of $14 \times 14 \times 3.5 \ \text{cm}^3$ so that it simulates the size of the top surface of the skull.

\noindent
\textit{Data generation:}
This grid was made up of $160 \times 160 \times 40$ grid points with a spatial separation of $8.75 \times 10^{-2} \ \text{cm}$ along all Cartesian coordinates.                                                                                                                                                                                                  Each surface of this grid was enclosed by a PML with 20 grid points, and an attenuation coefficient with a maximum value of 2 nepers per grid point was tapered within the PML \cite{Tabei}.
The pressure field was measured by $62 \times 62$ point-wise detectors, which were placed equidistantly on the top surface of the grid. The skull was simulated so that its top and bottom surfaces are aligned by the third and tenth horizontal planes of the grid points from the top surface of the cube. This has provided a thickness of 6.1mm for the skull, as well as a distance of 1.75mm between the top surface of the skull and the detection plane.

\noindent
\textit{Image reconstruction:}
Here, an inverse crime for discretisation was avoided by using a grid with different size, made up of $128 \times 128 \times 32$ grid points which are positioned with a spatial separation of 1.1mm along all Cartesian coordinates. Proportional to a reduction in size of the computational grid, we reduced the thickness of the PML to 16 grid points. Because of using a coarser computational grid, the thickness of the skull had to be reduced to 5.5mm with the top and bottom edges aligning the third and eighth horizontal planes of the grid points, respectively from the top surface of the grid.

\subsubsection{Physical parameters}

\noindent
Figures \ref{F4a}, \ref{F4b} and \ref{F4c}, respectively show the maps associated with the mass density, and the propagation speed of compressional and shear waves. As shown in these figures, the physical parameters of the medium are simulated the same as the 2D phantom.

\subsubsection{Validation of adjoint}
We used the inner product test in \eqref{adj_main} in order to numerically evaluate the accuracy of the computed adjoint model. To do this, we used a randomly selected vector for $\hat{P}$, together with a randomly chosen initial pressure distribution $P_0$ supported in the cuboid region below the skull. This is the region below the 10th (resp. 12th) horizontal plane from the top surface of the grid for image reconstruction (resp. data generation). The mean relative difference between the left-hand and right-hand sides of \eqref{adj_main} among 10 attempts was $2.21 \times 10^{-4}$ and $3.47 \times 10^{-5}$ for the grids used for image reconstruction and data generation, respectively.

\subsubsection{Simulation setting}

\begin{figure}\centering
{\subfigure[]{\includegraphics[scale=0.15]{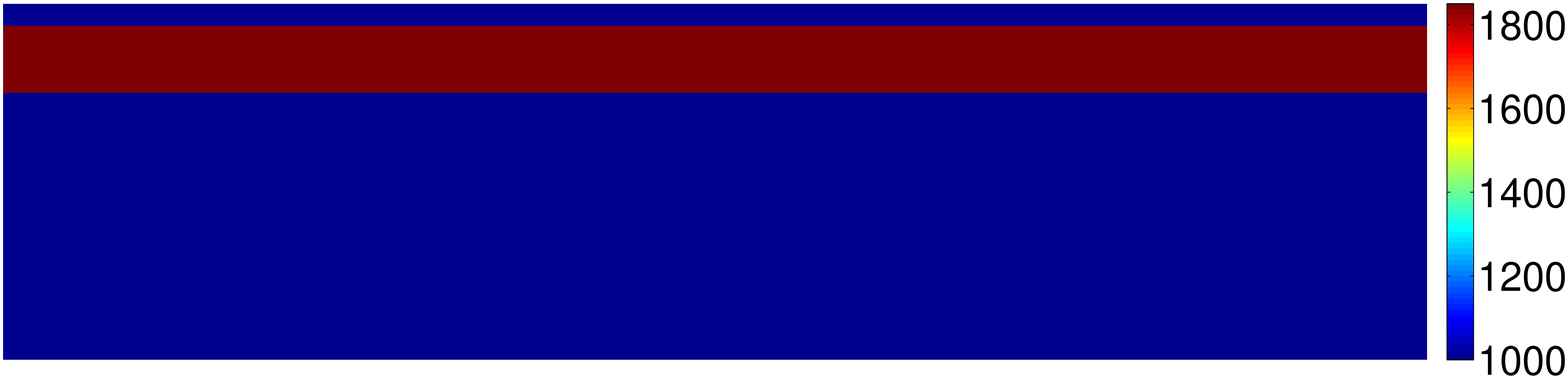}\label{F4a}}                                                                                                                                                       \hspace{0.25cm}                                                                                                                                                                                                                        \subfigure[]{\includegraphics[scale=0.15]{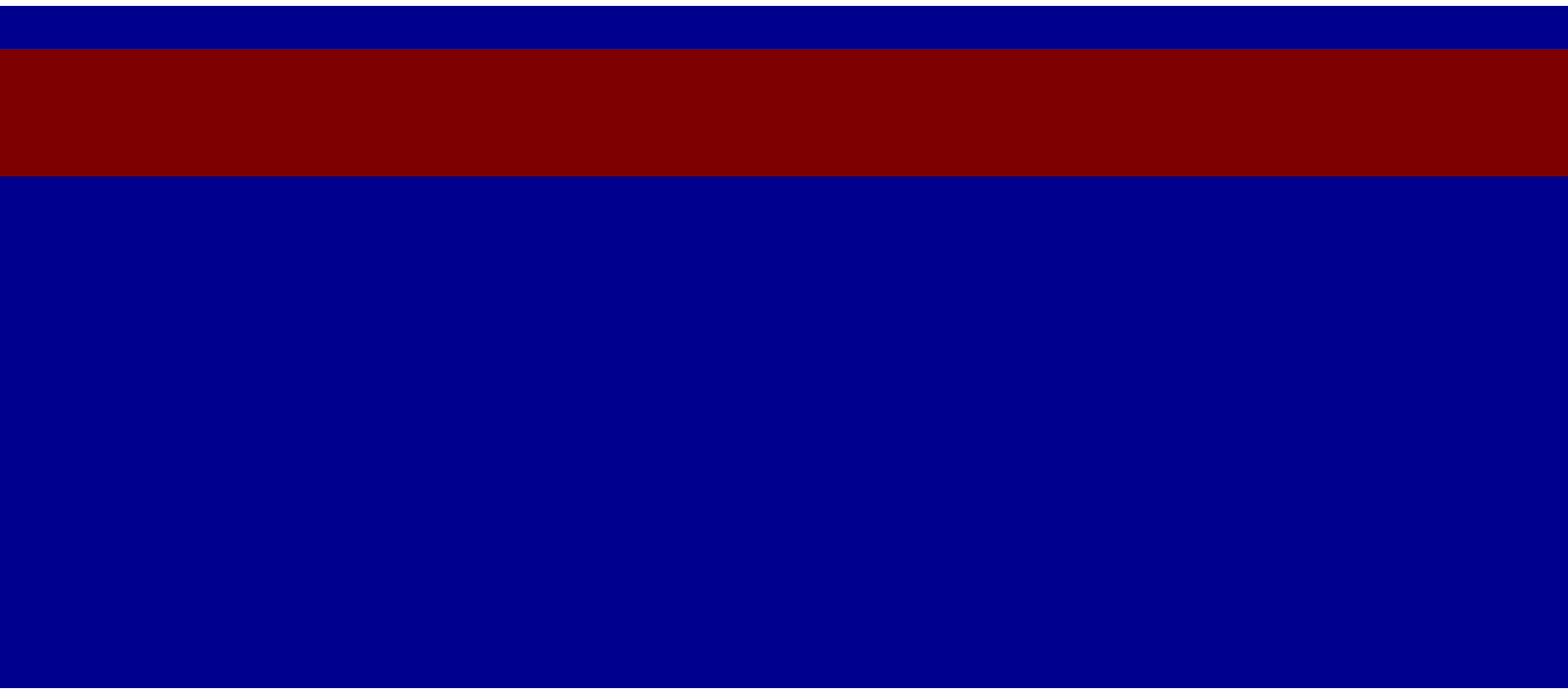}\label{F4b}}                                                                                                                                                           \hspace{0.25cm}                                                                                                                                                                                                                        \subfigure[]{\includegraphics[scale=0.15]{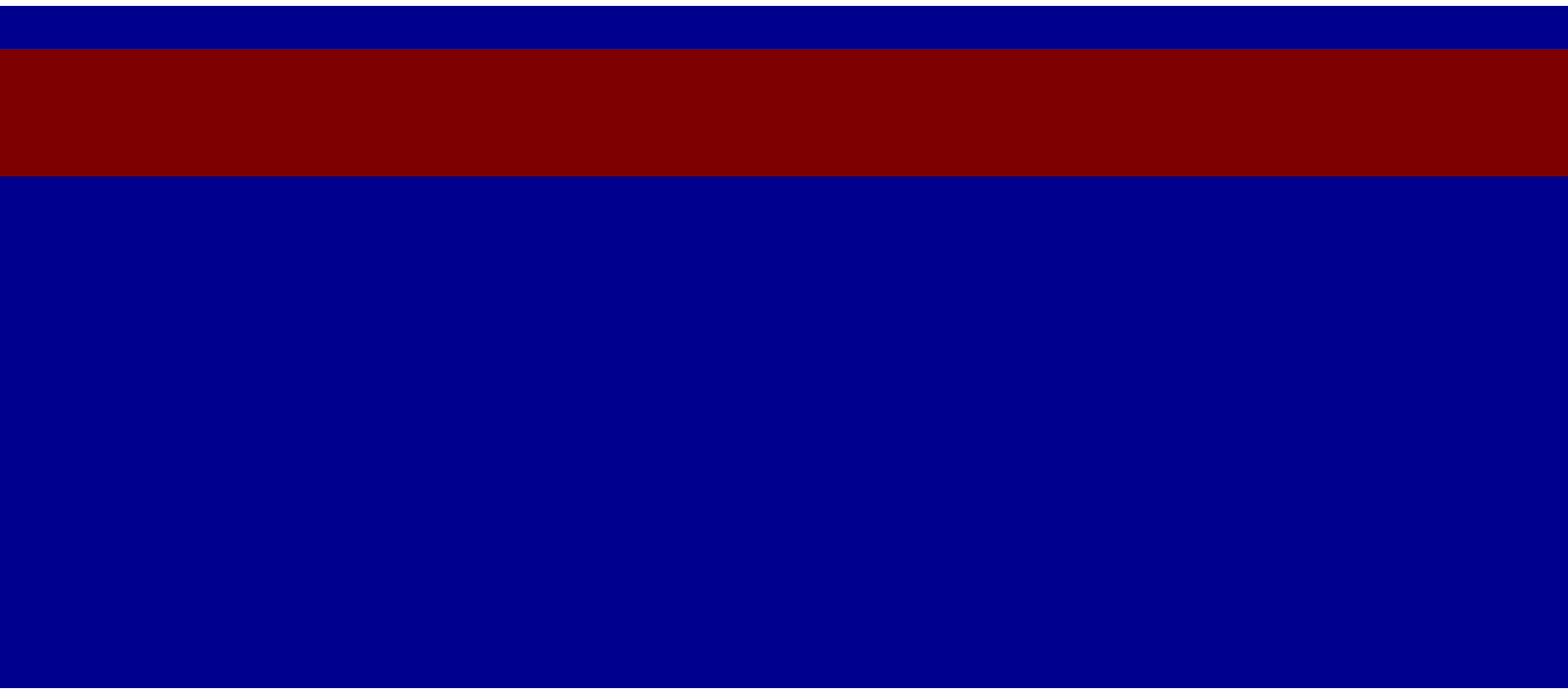}\label{F4c}}
}
\caption{3D phantom. Exact physical maps: (a) $\rho$ (b) $c_p$ (b) $c_s$.}
\end{figure}

\noindent
For image reconstruction, two scenarios were considered:

\noindent
\textit{Scenario 1:}
The maps that are displayed in figures \ref{F4a}, \ref{F4b} and \ref{F4c} were used for both data generation and image reconstruction. As discussed above, this is an inverse crime in estimating physical parameters, although the shift of soft tissue-skull interfaces between the fine and coarse grids cannot be neglected because of the high contrast between physical properties of the skull and soft tissue. Using these maps, the grid used for data generation supports a maximal frequency of $8.7514 \times 10^5 \text{Hz}$ for compressional waves across the entire medium and for shear waves propagated through the skull.

\textit{Scenario 2:}
In addition to the shifting error in physical parameters because of the discretisation,
these maps have been contaminated with a 30dB AWGN noise for data generation, whereas the reconstruction is done using the clean maps. This induces an error in estimating physical parameters, as they are not available exactly for image reconstruction. Using the noise contaminated maps, the grid used for data generation supports maximal frequencies up to $7.2018 \times 10^5 \text{Hz}$ and $7.4143 \times 10^5 \text{Hz}$ for compressional waves across the entire grid and shear waves within the skull, respectively.

The grid used for image reconstruction supports a maximal frequency of $6.8571 \times 10^5 \text{Hz}$ for compressional waves through the entire medium and shear waves propagated through the skull.
For simulating the initial pressure map, the phantom used for the 2D scenario was placed obliquely inside the cubic grid in a way in which the initial pressure distribution associated with the phantom is compactly supported in the soft tissue. Figure \ref{F4a} shows the simulated phantom from a top view. Using a $CFL$ of 0.3, the simulated pressure wavefield was recorded in $1532$ time steps. The recorded pressure field was then interpolated to the detectors using trilinear interpolation \cite{Mitsuhashi}. Similar to the 2D phantom, the vector of generated data $\hat{P}$ was contaminated with a $30$dB AWGN.
\subsubsection{Image reconstruction}
\begin{figure}\centering                                                                                                                                                                                                               {\subfigure[]{\includegraphics[width=.45\textwidth]{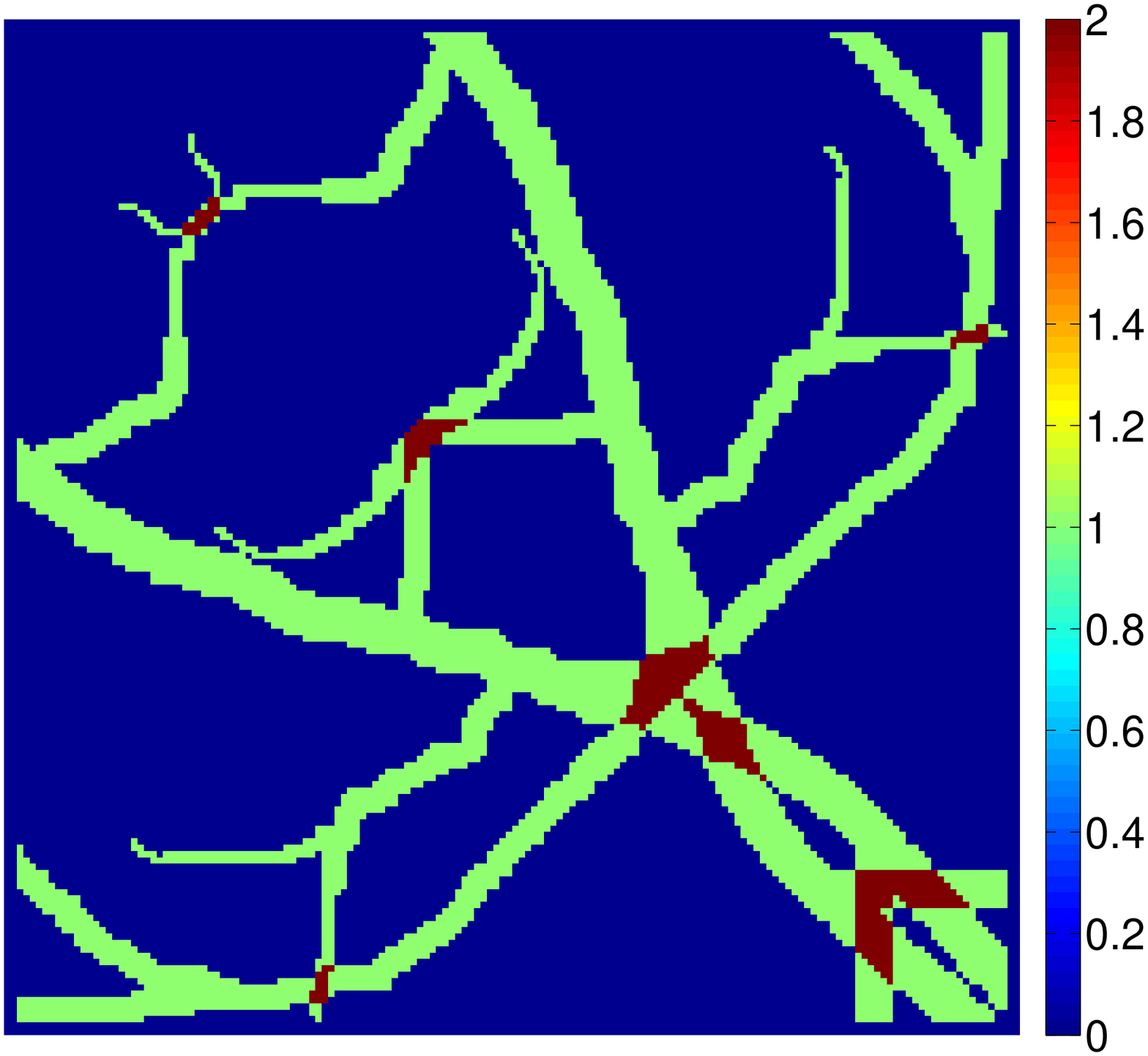}\label{F5a}}\\
\subfigure[]{\includegraphics[width=.45\textwidth]{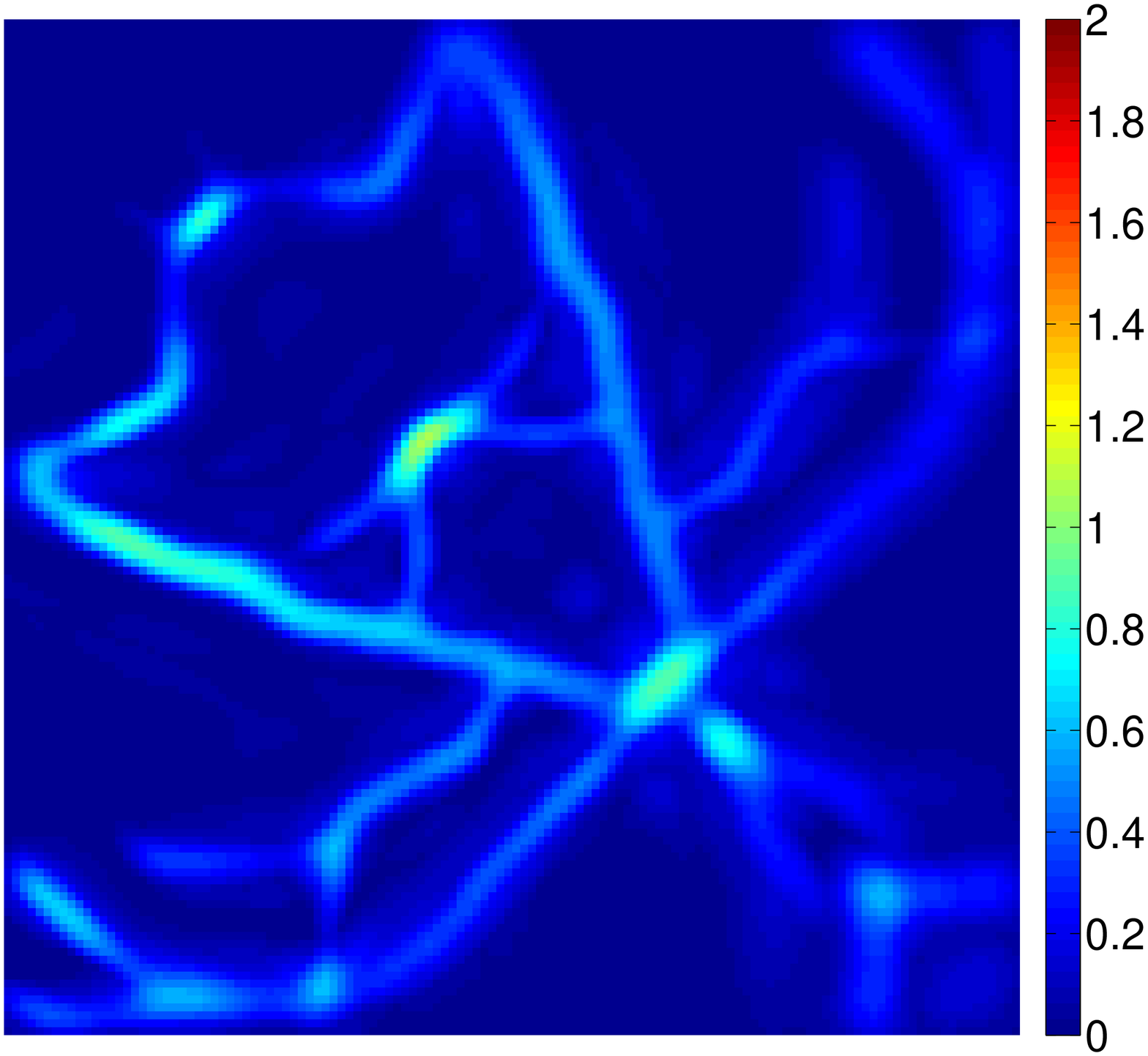}\label{F5b}}
\subfigure[]{\includegraphics[width=.45\textwidth]{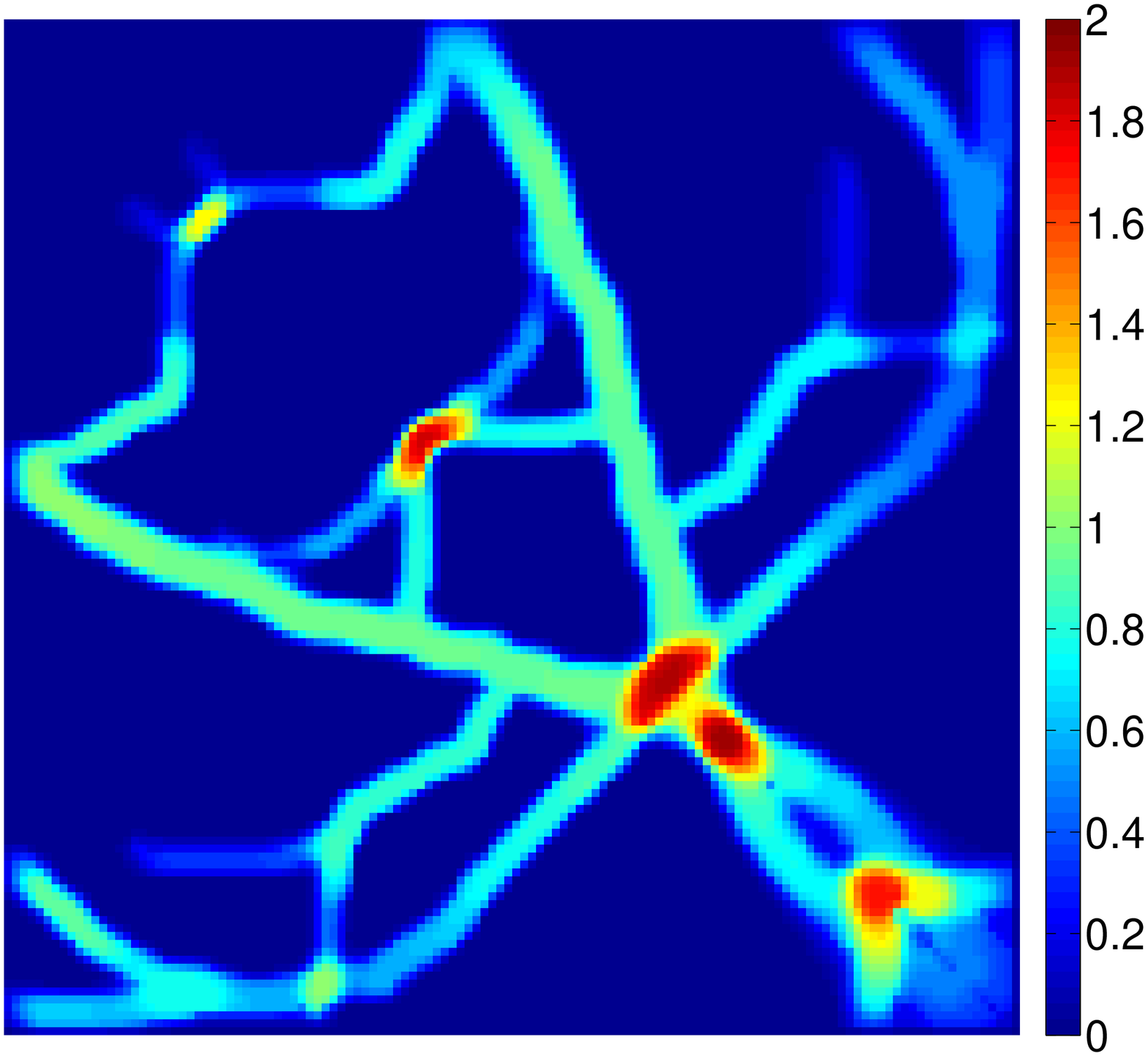}\label{F5c}}\\
\subfigure[]{\includegraphics[width=.45\textwidth]{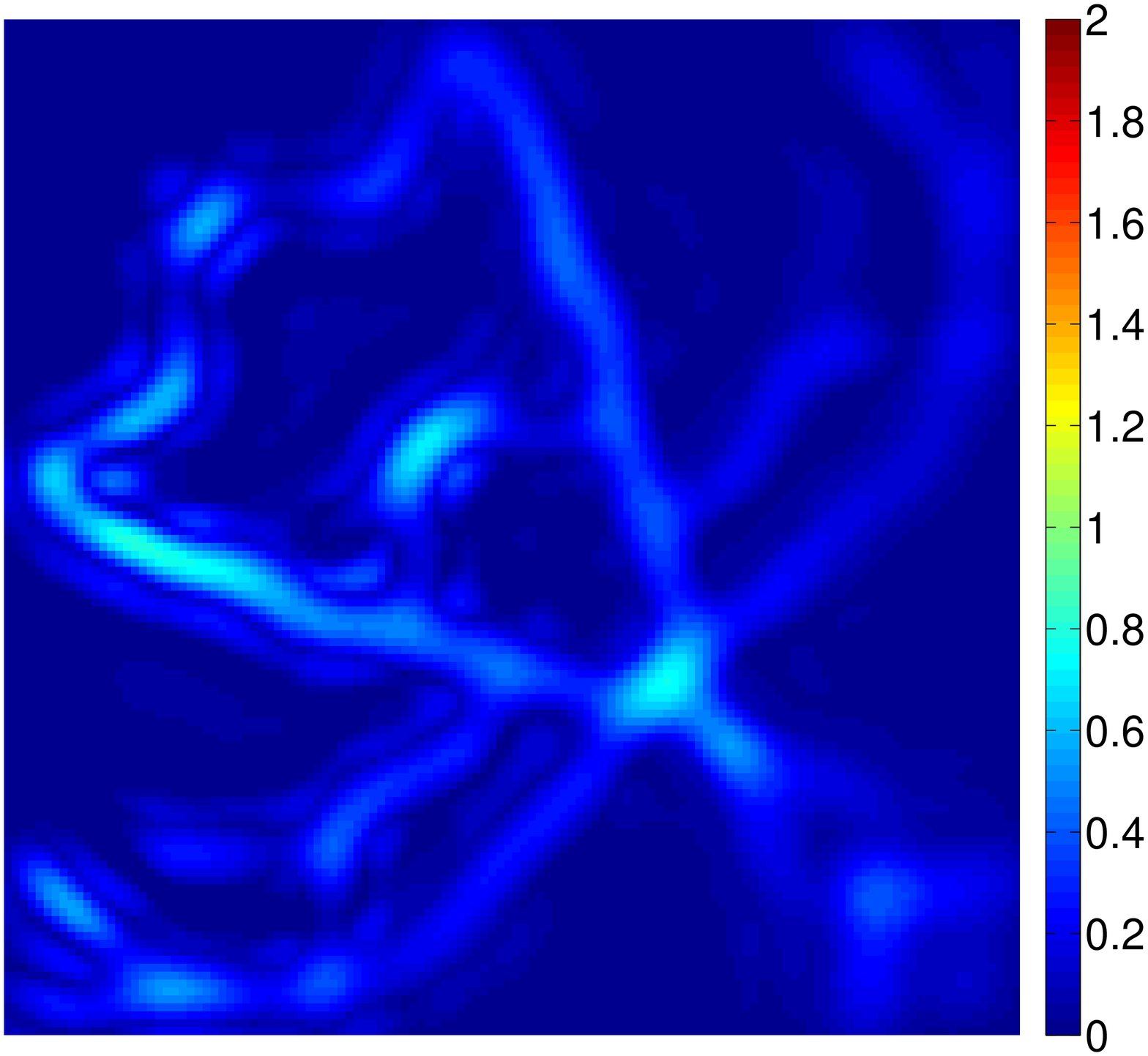}\label{F5d}}
\subfigure[]{\includegraphics[width=.45\textwidth]{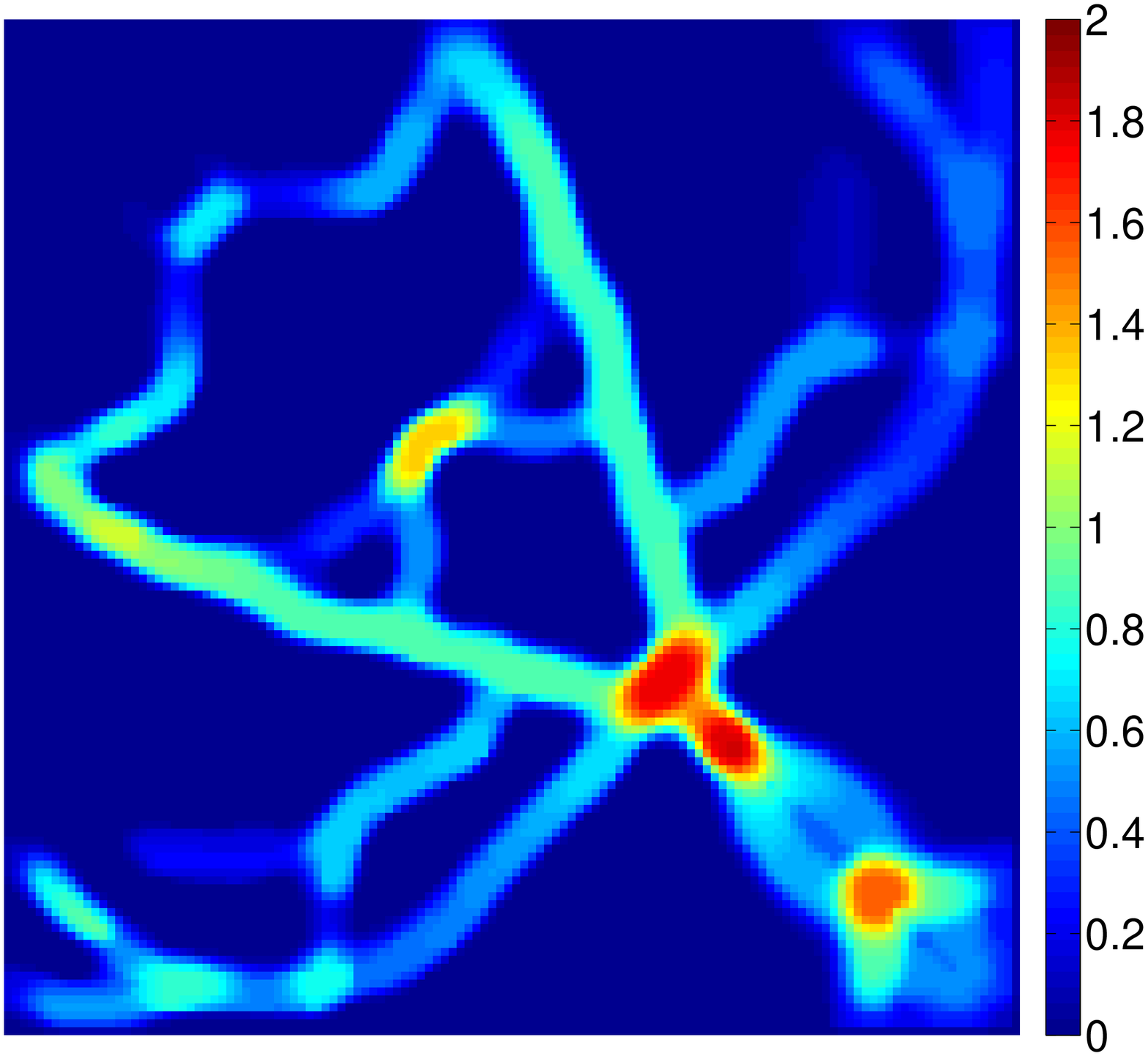}\label{F5e}}
}
\caption{3D phantom. (a) initial pressure map, and reconstructed images using exact physical parameters (inverse crime): (b) TR ($\alpha_{0,p,s}=0$) (c) ISTA, and erroneous physical parameters: (d) TR ($\alpha_{0,p,s}=0$) (e) ISTA.}
\end{figure}

We first reconstructed an image for each scenario using TR, which is available on the \textit{k-Wave} website \cite{treeby-b,Treeby-tool}. Here, all 3D images including phantom are displayed from a top view using \textit{maximum intensity projection} technique. The images reconstructed by TR for scenarios 1 and 2 are shown in figures \ref{F5b} and \ref{F5d}, respectively.
Using ISTA, the reconstruction parameters were chosen the same as for the 2D case.
We used the power iteration method for computing $L_f$. Figures \ref{F5c} and \ref{F5e} show an image of the final iterate computed by ISTA for scenarios 1 and 2, respectively.
A comparison between these two images indicates that using erroneous physical maps have led to a slight blurriness in the reconstructed image. Note that here the inverse crime has been avoided by exaggeration compared to real cases.

\subsubsection{Observations}

\begin{figure}\centering                                                                                                                                                                                                               {\subfigure[]{\includegraphics[width=.45\textwidth]{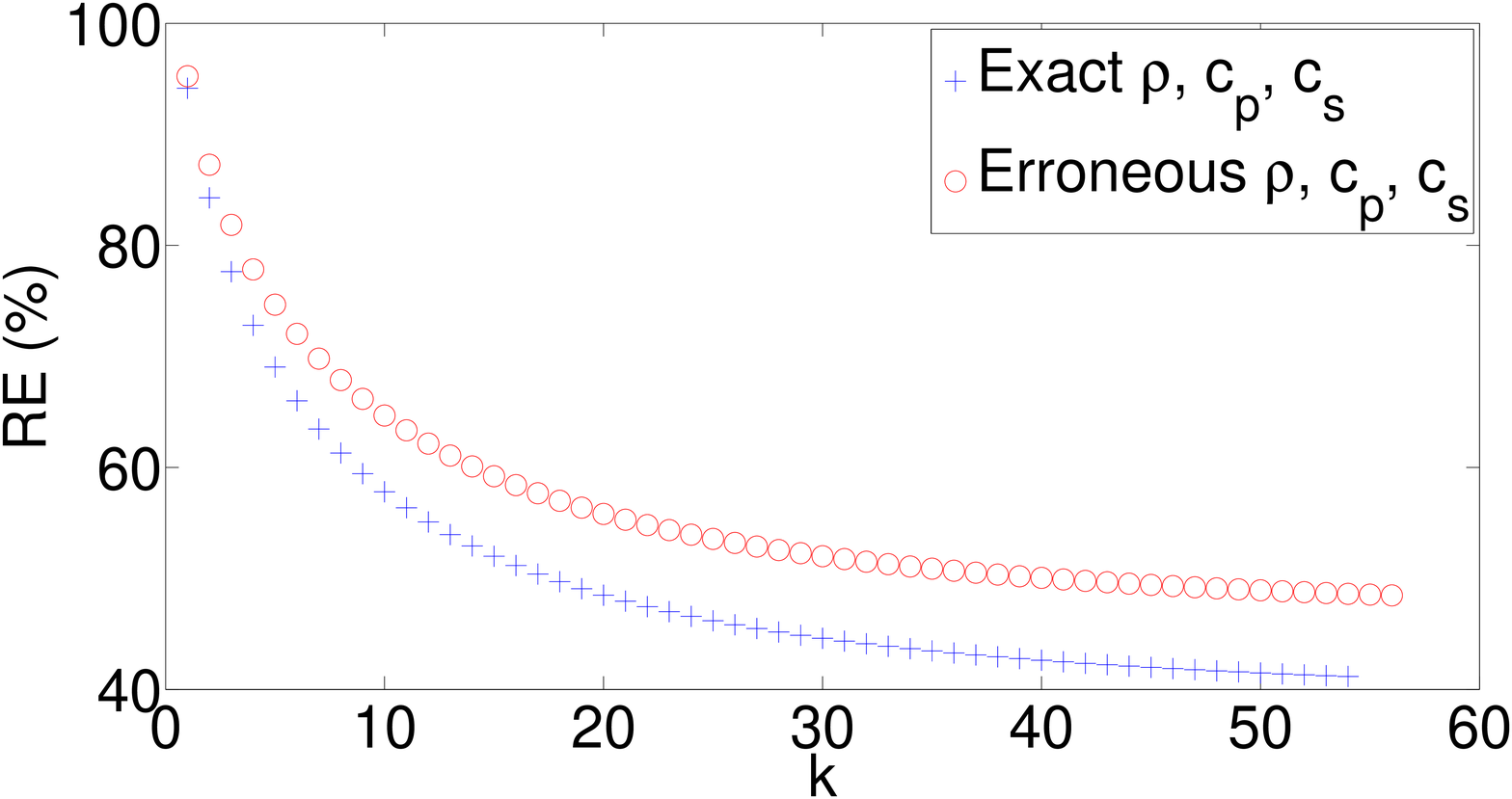}\label{F6a}}
\subfigure[]{\includegraphics[width=.45\textwidth]{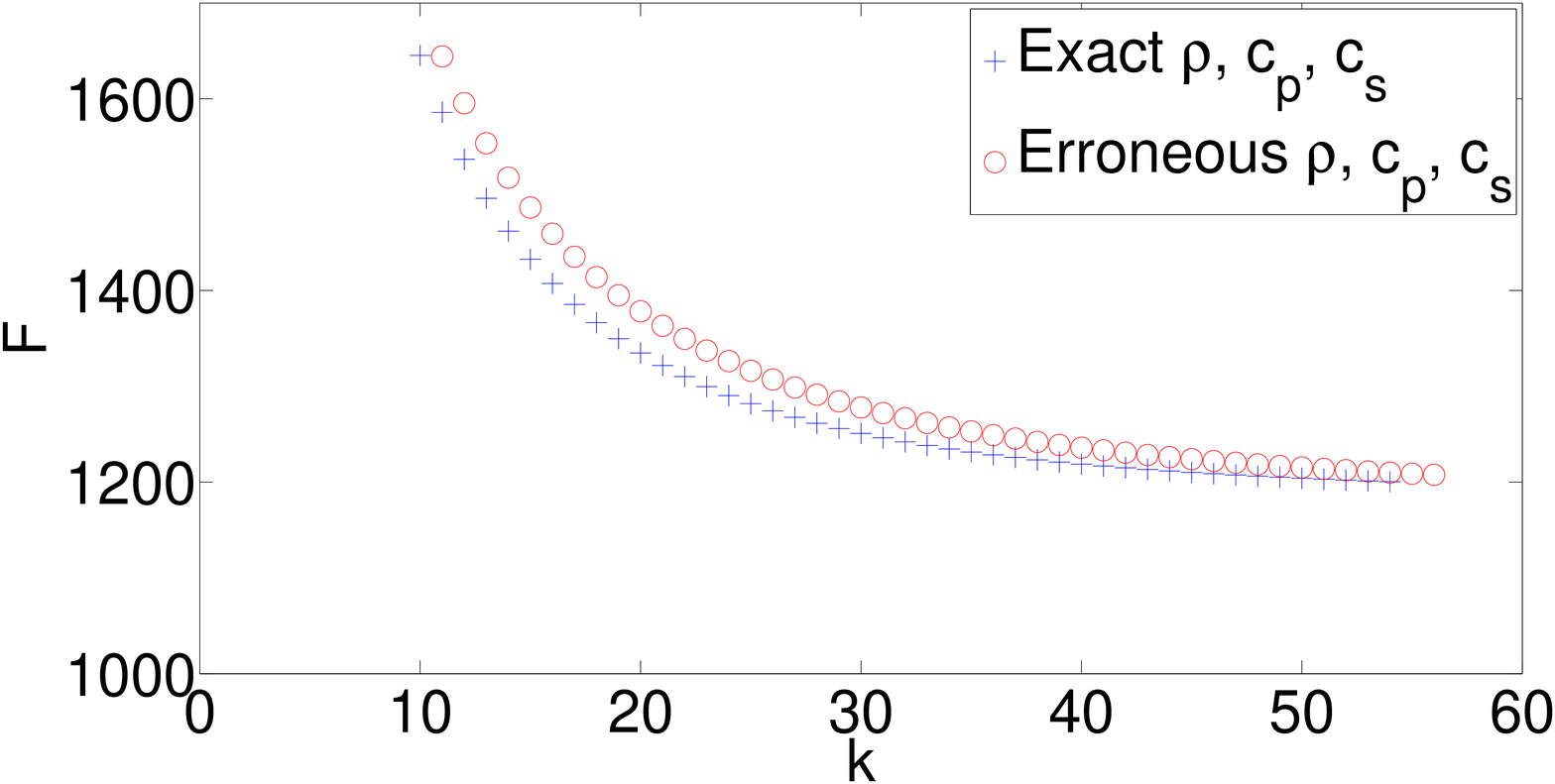}\label{F6b}}
}
\caption{3D phantom. (a) Relative Error (RE) (b) a large view of objective function F around the stopping point.}                                                                                                                                                                                                             \end{figure}

Figure \ref{F6a} shows the RE of the computed iterates versus iteration number. Additionally, figure \ref{F6b} shows the objective function values versus the iteration number around the terminating point. These plots have been displayed using the same colours as for the 2D phantom. From these, our observations for the two discussed scenarios are as follows.

\noindent
\textit{Scenario 1:}
Both $RE$ and $F$ monotonically decreased until the iteration $54$ at which the stopping criterion was satisfied. The final iterate, which is shown in figure \ref{F5c}, has an RE of $41.13\%$ and an $F$ of $1.19 \times 10^3$.

\noindent
\textit{Scenario 2:}
Using the noise contaminated physical maps for data generation, a monotonic reduction for $RE$ and $F$ is observed again, and the terminating criterion was satisfied at iteration $56$. The final reconstructed image, which is shown in figure \ref{F5e}, has an RE of $48.44\%$ and an $F$ of $1.20 \times 10^3$.

\section{Discussion and conclusion}
In this work, we derived the adjoint of the continuous map defined in \eqref{elast_map} and \eqref{init1}, which describes the propagation of PA waves in linear isotropic viscoelastic media with the absorption and physical dispersion following a frequency power law. We analytically showed that a numerical computation of our continuous adjoint using a k-space pseudo-spectral method matches the algebraic adjoint of an associated discretised map defined by \eqref{dip} and \eqref{for_map_dis}.

From a numerical point of view, it was shown that this forward and adjoint pair satisfies the inner product test in \eqref{adj_main}. This pair was then incorporated in a positivity constrained optimization algorithm based on ISTA that is regularized by the TV denoising approach of Chambolle \cite{Beck}. We preferred to test the derived forward and adjoint operators on a classical inverse solver (cf. \cite{Arridge,Javaherian} for the application of ISTA in PAT), although this poses some limitations such as a low speed of convergence. The convergence can be improved by using a a fast version of ISTA (FISTA) \cite{Huang,Arridge-b,Javaherian}. In addition, an iterative enhancement of solutions based on a \textit{Bregman iteration} \cite{Osher} may be useful when compressively sampled data are collected. A \textit{Bregman iteration} algorithm using FISTA has been successfully applied in this case \cite{Arridge-b}.

In both 2D and 3D cases, the iterates are monotonically converged to a minimizer of an objective function, and the final reconstructed image was close to the ground truth image. In presence of some levels of error in estimating physical parameters, the iterates are monotonically converged again, but the iterate at the stopping point was slightly less accurate than using the exact physical parameters. This loss of accuracy cannot be detected by eye in the 2D scenario, as shown in figure \ref{F2e}. However,
figure \ref{F5e} shows that an error in estimating physical parameters has led to a slight loss of contrast in the 3D scenario, compared to using exact physical maps. Note that in the 3D scenario for the grid used for data generation, in addition to a 30 dB noise added to the physical maps, the skull's thickness is 0.6mm larger than the grid for image reconstruction.

In addition, for the 2D scenario, as shown in figures \ref{F2c} and \ref{F2e}, the reconstructed images have some blurriness in regions close to the skull. We believe that this can be attributed to the full internal reflection of wavefronts nearly tangent to the skull, and agrees with theoretical predictions of stability for inversion found in \cite{Stefanov} using methods of microlocal analysis. In essence, the blurred region close to the skull is not fully resolved because the wavefronts emanating from that region do not reach the detectors (Note that the geometry of skull in our study is not realistic). To fully understand this a more delicate study on the relation between the theoretical analysis of \cite{Stefanov} and PAT of the brain using optimization algorithms may be needed.

The simplified geometries of the skull we used in our simulations look sufficient to provide an insight on the performance of the derived adjoint, but the geometry and composition of the skull in real cases are much more complicated than our simulations \cite{Huang-a}. In addition, in practical cases, to extract all information available from the measured data, the maximal frequency supported by the computational grid must match the maximal frequency that is detectable by detectors \cite{Treeby}. This dramatically increases the computational demands regarding storage space and speed, but it can be handled using GPU accelerated computing \cite{Kang}, or Field-programmable gate array (FPGA) \cite{Alqasemi}. The 3D detection setting in our study simulated a planar Fabry-P\'erot (FP) photoacoustic scanner, which requires several minutes to collect time series of data from PA wavefields \cite{Arridge-b}. Further studies can be done to apply our optimization algorithm on ultra-fast PAT acquisition systems that utilize spatio-temporal sub-sampled data \cite{Arridge-b}.

Using our derived adjoint, an extension of the PAT problem of brain to direct quantitative PAT (QPAT), a direct estimation of the optical parameters inside the skull from the acoustic data collected outside the skull, would be a very interesting topic. The arising \textit{opto-elastic} inverse problem is more challenging than the \textit{opto-acoustic} problem because of the high optical absorption and scattering of the skull and low degrees of freedom for optical illumination. This limits the applicability of multi-source QPAT, which is necessary for uniqueness of the problem when we use a single-frequency optical excitation \cite{Haltmeier}.

\section{Appendix}
In this appendix, we will show that the adjoint operator can be put into the form of a system of coupled partial differential equations, in the same way as the forward operator, and the update of
particle velocity field is actually a sum of the adjoint of absorption and dispersion terms enforced to the stress tensor field.

\textit{Continuous adjoint:}
This is derived by plugging the first formula in \eqref{stress2:eq} into \eqref{vel2:eq}
\begin{align} \label{vell:eq}
\begin{split}
 & \rho \frac{{{\partial \mathbf{v}_i}^{p,s}}^*}{\partial t}=  \sum_{p,s} q^{p,s}  \Bigg( \frac{\partial}{\partial x_i} \Big (\lambda {\sigma_{ll}^{p,s}}^* \Big ) + 2 \frac{\partial}{\partial x_j} \Big( \mu {\sigma_{ij}^{p,s}}^* \Big)  \Bigg)+\\
&  \Big( \sum_{p,s} q^{p,s} \sin(\pi y/2) L_{c_{p,s}}^{(y-1)*}-\cos(\pi y/2) L_{c_{p,s}}^{(y-2)*} \frac{\partial}{\partial t}  \Big)   \Bigg( \frac{\partial}{\partial x_i} \Big
(\chi {\sigma_{ll}^{p,s}}^* \Big ) + 2 \frac{\partial}{\partial x_j}  \Big( \eta {\sigma_{ij}^{p,s}}^* \Big)\Bigg)
\end{split}
\end{align}

\textit{Discretised adjoint:}
In \eqref{main_dis2}, plugging the second line into the first line yields
\begin{align}  \label{vell:dis}
\begin{split}
{\tilde{\mathbf{v}}}_{n+1/2}  & = A_v  \left( A_v \tilde{\mathbf{v}}_{n-1/2}+\Psi_\text{dis}^* \tilde{\sigma}_{n}-\frac{1}{\Delta t} \Big( \big(  \Psi_\text{abs}^{\prime^*}\tilde{\sigma}_n \big)^{p,s} -A_v^2  \big( \Psi_\text{abs}^{\prime^*}\tilde{\sigma}_{n-1}\big)^{p,s} \Big) \right).
\end{split}
\end{align}

The numerical computation of \eqref{vell:eq} is the same as \eqref{vell:dis}, except how the PML acts on the temporal gradient of the stress tensor field. These formulae require an explicit computation of the temporal gradient of the stress tensor using finite difference schemes. To avoid this, we used the formulae \eqref{vel2:eq} and \eqref{stress2:eq} (resp. \eqref{main_dis2}) for the continuous (resp. discretised) adjoint, which are computed the same, as discussed in sections \ref{nmethodf} and \ref{s6}, respectively.

\bibliographystyle{abbrv}                                                                                                                                                                                                              \bibliography{my_ref}

\end{document}